\documentclass[11pt,a4paper,twoside]{m-art}

\usepackage{times}
\usepackage{a4}
\usepackage{xy} 
\xyoption{all}
\usepackage{amsmath, amssymb}

\def\Char{\mathrm{char}}
\def\CH{\mathrm{CH}}
\def\ev{\mathrm{ev}}

\def\id{\mathrm{id}}
\def\Mor{\text{Mor}}
\def\Spec{\mathrm{Spec}}
\def\Proj{\mathrm{Proj}}
\def\virt{\mathrm{virt}}
\def\Isom{\underline{\mathrm{Isom}}}

\def\abs#1{\lvert#1\rvert}

\swapnumbers
\newtheorem{Def-s}[subsection]{Definition}
\newtheorem{Prop-s}[subsection]{Proposition}
\newtheorem{Thm}[subsubsection]{Theorem}
\newtheorem{Def}[subsubsection]{Definition}
\newtheorem{Rem}[subsubsection]{Remark}
\newtheorem{Prop}[subsubsection]{Proposition}
\newtheorem{PropDef}[subsubsection]{Proposition and Definition}
\newtheorem{Lem}[subsubsection]{Lemma}
\newtheorem{Cor}[subsubsection]{Corollary}

\newtheorem{DefRem}[subsubsection]{Definition and Remark}

\def\C{\mathbb{C}}

\def\P{\mathbb{P}}
\def\Q{\mathbb{Q}}
\def\R{\mathbb{R}}
\def\Z{\mathbb{Z}}

\def\AA{\mathcal{A}}
\def\BB{\mathcal{B}}
\def\CC{\mathcal{C}}
\def\DD{\mathcal{D}}
\def\LL{\mathcal{L}}

\def\OO{\mathcal{O}}
\def\PP{\mathcal{P}}
\def\RR{\mathcal{R}}
\def\WW{\mathcal{W}}

\def\AAA{\mathfrak A}

\def\GGG{\mathfrak G}

\def\VVV{\mathfrak V}

\def\Mbar{\overline{\mathcal{M}}}

\newcommand\stv[2]{\left.\kern-\nulldelimiterspace
                \left\{#1\vphantom{#2}\,\right|#2\right\}}

\begin{document}

\author{Arend Bayer}
\address{Arend Bayer, Max-Planck-Institut f\"ur Mathematik, Bonn, Germany}
\email{bayer@mpim-bonn.mpg.de}
\curraddr{Mathematical Sciences Research Institute, 17 Gauss Way, Berkeley,
CA 94708}
\author{Yu.~I.~Manin}
\address{Yuri I.~Manin, Max-Planck-Institut f\"ur Mathematik, 
Vivatsgasse 7, Bonn, Germany; and
Northwestern University, Evanston, USA}
\email{manin@mpim-bonn.mpg.de}

\title{Stability Conditions, Wall-crossing 
and weighted Gromov-Witten Invariants}


\begin{abstract}
We extend B.~Hassett's theory of  weighted
stable pointed curves (\cite{Has03}) to weighted stable maps.
The space of stability conditions is described
explicitly, and the wall-crossing phenomenon studied.
This can be considered as a non-linear analog of
the theory of stability conditions in abelian and triangulated categories
(cf. \cite{Christmas}, \cite{Bridgeland:Stab}, \cite{Joyce1, Joyce2,
Joyce3, Joyce4}).

We introduce virtual fundamental classes and thus obtain weighted
Gromov-Witten invariants. We show that by including gravitational
descendants, one obtains an $\LL$-algebra as introduced in
\cite{LosevManin-Ext} as a
generalization of a cohomological field theory.
\end{abstract}

\dedicatory{\`a Pierre Deligne, en t\'emoignage d'admiration}
 
\maketitle

\setcounter{section}{-1}

\section{Introduction: Hassett's stability conditions}

\subsection{Pointed curves} {\it A nodal curve $C$} over an algebraically
closed field $k$ is a proper nodal reduced one-dimensional 
scheme of finite type over this field whose only singularities are nodes.
The genus of $C$ is $g:=\dim\,H^1(C,\mathcal{O}_C).$

Let $S$ be a finite set. {\it A nodal $S$-pointed curve $C$}
is a system $(C, s_i\,|\,i\in S)$ where $\{s_i\}$ is a family of closed
non-singular $k$-points of $C$, not necessarily pairwise distinct.
The element $i$ is called the label of $s_i$.

The normalization $\widetilde{C}$ of $C$ is a disjoint union
of smooth proper curves. Each irreducible component of
$\widetilde{C}$ carries inverse images of some labeled points $s_i$
and of singular points of $C$. Taken together, these points are called
{\it special ones.} Instead of passing to the normalization,
we may consider branches (local irreducible germs) of $C$ passing
through labeled or singular points. They are in a natural bijection
with special points.

A  nodal connected $S$-pointed curve 
$(C, s_i)$ is called {\it stable} if $s_i\ne s_j$
for $i\ne j$ and any of the 
following three equivalent 
conditions hold:

(i) {\it The automorphism group of $(C, s_i)$ is finite.}

(ii) {\it Each irreducible component of
$\widetilde{C}$ of genus 0 (resp. 1) supports $\ge 3$
(resp.  $\ge 1$) distinct special points.}

(iii) {\it The line bundle 
$\omega_C\left(\sum_{i\in S} s_i\right)$ is ample.}

This definition has a straightforward extension to families
of stable $S$-pointed curves (cf. below). The basic result
states that families of stable
$S$-pointed curves of genus $g$ form (schematic points of) a connected
smooth proper over $\mathbb{Z}$
Deligne-Mumford stack $\overline{\mathcal{M}}_{g,S}.$ 
It contains an open dense substack
${\mathcal{M}}_{g,S}$ parameterizing irreducible smooth curves,
and is its compactification.

\subsection{Weighted stability.} Generalizing condition (iii),
B.~Hassett enriched
the theory by additional parameters generating a whole new family
of stability conditions, which lead to new moduli
stacks, representing different compactifications of ${\mathcal{M}}_{g,S}$.

Namely,  {\it the weight data} on $S$ is
a function $\mathcal{A}:\,S\to \mathbb{Q}$, $0 < \AA(i)\le 1.$
$S$ together with a weight data will be called {\it a weighted set.}

\begin{Def}[\cite{Has03}] 			\label{def:weighted-stable}
A connected $S$-pointed curve 
$(C, s_i\,|\,i\in S)$ is called 
weighted stable (with respect to $\mathcal{A}$)
if the following conditions 
are satisfied:

(i) $K_{C} + \sum_i \AA(i)s_i$ is an ample divisor, where $K_{C}$ is the 
canonical class of $C$.

(ii) For any subset $I\subset S$ such that $s_i$ pairwise coincide
for $i\in I$, we have $\sum_{i\in I} \AA(i) \le 1.$
\end{Def}

Clearly, (i) implies that $2g-2+\sum_i \AA(i) >0.$ 

The usual stability
notion corresponds to the case $\AA(i)=1$ for all $i\in S$.
Independently of Hassett's work, A.~Losev and Yu.~Manin
considered in \cite{LosevManin00}, \cite{LosevManin-Ext} some
non-standard moduli spaces which turned out to correspond
to special Hassett's stability conditions: see \cite[section 6.4]{Has03}
and \cite{Manin-LgS}.  

Definition \ref{def:weighted-stable} admits a straightforward extension
to families:

Let $U$ be a scheme, $S$ a finite set, $g\ge 0$. 
{\it An $S$-pointed nodal  curve (or family of curves) of genus $g$
over $U$} consists of the data
$$
(\pi:\,C\to U;\,s_i:\,U\to C,\ i\in S)
$$
where $\pi$ is a flat proper morphism whose geometric fibers $C_t$
are nodal $S$-pointed curves of genus $g$.

This family is called $\mathcal{A}$-stable iff

(i) {\it $K_\pi + \sum_i \AA(i)s_i$ is $\pi$-relatively ample.}

(ii) {\it For any $I\subset S$ such that 
$\cap_{i\in I} s_i\ne \emptyset$,  we have $\sum_{i\in I} \AA(i) \le 1.$}

\subsection{Stacks of weighted stable curves
$\overline{\mathcal{M}}_{g,\mathcal{A}}$} 
The first main result of \cite{Has03}
is a proof of the following fact. Fix a weighted set of labels $S$
and a value of genus $g$. Then families of weighted stable
$S$-pointed curves of genus $g$ form (schematic points of) a connected
smooth proper over $\mathbb{Z}$
Deligne-Mumford stack $\overline{\mathcal{M}}_{g,\mathcal{A}}.$ The respective
coarse moduli scheme is projective over $\mathbb{Z}$.

\subsection{Walls and wall-crossing.} The further results of 
Hassett on which we focus in this introduction
concern the geometry of the space of stability 
conditions governing the varying geometry of boundaries of
$\overline{\mathcal{M}}_{g,\mathcal{A}}$ (\cite{Has03}, sec. 5).

Put
$$
D_{g,S}:=\{\mathcal{A} \in \mathbb{R}^S\,|\, 0<\AA(i)\le 1,\ \sum_{s} \AA(i) > 2-2g\}.
$$
{\it Walls} are non-empty intersections of $D_{g,n}$ with certain 
hyperplanes indexed by subsets $I\subset S$: 
$$
w_I:=\{\mathcal{A} \in D_{g,S}\,|\, \sum_{i\in I} \AA(i)=1\,\}.
$$
{\it Coarse chambers} are defined as connected components of
$$
D_{g,S} - \bigcup_{2< |I| \le n-3 \delta_{g, 0}} w_I .
$$
{\it Fine chambers} are connected components of
$$
D_{g,S} - \bigcup_{2\le |I|\le n-2 \delta_{g, 0}} w_I .
$$

B.~Hassett proves the following result:

\begin{Prop} (i) The moduli stack
$\overline{\mathcal{M}}_{g,\mathcal{A}}$
is constant on each coarse chamber, and differs from one
coarse chamber to another.

(ii) The universal curve $\overline{\mathcal{C}}_{g,\mathcal{A}}$
is constant on each fine chamber, and differs from one
fine chamber to another.

Finally, for any point $\mathcal{A}'$ belonging to a wall, there exists a point
$\mathcal{A}$ inside a neighboring coarse (resp. fine) chamber at which
$\overline{\mathcal{M}}_{g,\mathcal{A}}$ (resp.
$\overline{\mathcal{C}}_{g,\mathcal{A}}$)
is the same as at $\mathcal{A}'$.
\end{Prop}

\subsection{Plan of this paper.} Let $V$ be a smooth projective
manifold. M.~Kontsevich has defined stacks $\mathcal{M}_{g,S}(V)$
of $S$-pointed stable maps $(C \to V; s_i)$. The stacks
 $\mathcal{M}_{g,S}$ correspond to the case $V =$ a point.
In this paper we generalize Hassett's stability conditions
to $\mathcal{M}_{g,S}(V)$ and study the resulting stacks.

In \ref{geometry}, we define the precise moduli problem and 
construct its moduli space as a proper Deligne-Mumford stack. 
We show the existence of birational contraction morphism for any
reduction of the weights; in particular, all moduli spaces of
weighted stable maps are birational contractions of the
Kontsevich moduli space.

We establish the existence of 
all basic morphisms (gluing, changing the target, forgetting sections etc.)
between them in \ref{sect:elementary-morphisms}. Section
\ref{sect:chamber-decomposition}
describes the chamber decompositions of the set of admissible weights.
and exhibits the reduction morphisms for a wall-crossing as an explicit
blow-up.

In \ref{sect:virtual}, we postulate a list of basic
properties for virtual fundamental classes, and discuss consequences
for the weighted Gromov-Witten invariants.
After introducing the language of weighted graphs in
\ref{sect:graph-language}, we prove a more complete
graph-level list of properties of the virtual fundamental classes in
\ref{sect:virtual-and-graphs}. 

One motivation of this study was the work by Losev and Manin
on painted stable curves \cite{LosevManin00, LosevManin-Ext, Manin-LgS},
which constitute a special case of weighted stable curves.
The authors introduced the
notion of an $\LL$-algebra as an extension of the notion of a
cohomological field theory of \cite{KoMa}.

The construction of
virtual fundamental classes in the extended context
of new stability conditions allows us to produce Gromov-Witten invariants
based on weighted stable maps. Including gravitational descendants,
we obtain $\LL$-algebras in the sense of \cite{LosevManin-Ext}.
While weighted Gromov-Witten invariants without gravitational descendants
yield nothing new (see proposition \ref{prop:GW-without-gravity}),
the coupling to gravity in the weighted case
exhibits a new structure on quantum cohomology.

In \cite{Mustatas:Chowring}, the authors already constructed moduli spaces of
weighted stable maps as an aide in computing the Chow ring of non-weighted
stable maps with target $\P^n$.
Independently of the present paper, Alexeev and Guy studied the behavior
of gravitational descendants for changes of weights in \cite{Alexeev-Guy},
assuming the same definition of virtual fundamental classes that we study 
in sections \ref{sect:virtual} and \ref{sect:virtual-and-graphs}.

Another motivation is spelled down below.

\subsection{Stability conditions in abelian and triangulated
categories.} Stability conditions have been generally designed to
choose a {\it preferred} compactification of various moduli
spaces, typically of vector bundles, or more general
coherent sheaves on projective manifolds. It was only recently
that the attention of algebraic geometers shifted
to the families of variable stability conditions and their geometry:
see \cite{Christmas}, \cite{Bridgeland:Stab}, \cite{Joyce1, Joyce2,
Joyce3, Joyce4}, and the references
therein.  An influential recent paper by T.~Bridgeland \cite{Bridgeland:Stab}
was very much stimulated by
physics work on mirror symmetry, in particular, M.~Douglas's
notion of $\Pi$-stability.

In this subsection we will sketch a purely geometric context
in which various notions of stability in derived categories
of coherent sheaves might be quite useful (see \cite{Inaba},
\cite{Bridgeland:Flop},
for a version of background notions). 
   
Namely, consider the problem treated in several papers
by A.~Bondal, D.~Orlov and others: {\it what can be said
about a (smooth projective) manifold $V$ if we know its
bounded derived category
of coherent sheaves $D(V)$?}

In an important paper \cite{Bondal-Orlov} it was shown that if
the canonical sheaf $\Omega_V$ of $V$ or its inverse is ample, then
$V$ can be reconstructed up to an isomorphism from
$D(V)$. The strategy of proof is this:
the authors show how to detect (up to a shift) 
classes of structure sheaves
of closed points of $V$ in $D(V)$, then classes
of invertible sheaves, and finally to reconstruct
the canonical (or anticanonical) homogeneous coordinate ring.

This result can become dramatically wrong, when $\Omega_V^{\pm 1}$ 
is not ample, for example, when it is trivial.
In the proper Calabi-Yau case various birational
models may lead to equivalent derived categories.
The complete picture in this case is far from being clear.
The proof that worked in the Fano/anti-Fano cases
breaks down at the first step: the classes
of structure sheaves
of closed points of $V$ become unrecognizable.

However, the general strategy of the proof could be saved
without additional assumptions on $\Omega_V$  if one
could do the following:

a) Devise a family of appropriate stability conditions $\mathcal{C}$
(this is probably already done in \cite{Bridgeland:Stab}).

b) Prove that various $V$'s with ``the same'' $D(V)$
could be reconstructed as moduli spaces $V_{\mathcal{C}}$ of appropriately
defined $\mathcal{C}$-stable point-like complexes
in $D(V)$. The deformation theory of objects in derived categories
is not yet a mature subject, but see \cite{Lunts-Orlov:DG-deformation}
for some recent developments.

c) obtain a sufficiently detailed description of
chamber decomposition and wall-crossing in the space of $\mathcal{C}$'s.

A tentative picture of this type can be glimpsed from
the Aspinwall's sketch \cite{Aspinwall:points}.
Locally, the wall-crossing phenomenon has been studied
in \cite{Toda:stab-crepant_res}.

Notice however that it is not clear a priori
what would be the net outcome of such a reasoning.
In fact, according to the recent preprint
\cite{Caldararu:Kawamata-counterexample},
two Calabi-Yau threefolds can have equivalent
derived categories without being birationally equivalent.

On the positive side, however, they must have
isomorphic motives: cf. \cite{Orlov:derived-and-motives}.

From this perspective, Hassett's theory and its generalization, 
discussed in this paper, can be perceived as a toy model
for the more sophisticated case of the triangulated categories.
Moreover, various notions of stability for maps of curves
into nontrivial target spaces could conceivably be combined
with similar stability notions for complexes of sheaves 
on the target space leading
to a richer structure of quantum cohomology.

\subsection{Acknowledgements}
The first author would like to thank Andrew Kresch for useful remarks on
virtual fundamental classes.

\section{Geometry of moduli spaces of weighted stable maps}
\label{geometry}

\subsection{The moduli problem}
Let $k$ be a field of any characteristic, $V/k$ a projective variety,
and $\beta \in \CH^1(V)$ an effective one-dimensional class in the
Chow ring.
Let $S$ be a finite set with weights $\AA \colon S \to \Q \cap [0, 1]$,
and let $g \ge 0$ be any genus. 

\begin{Def}					\label{def:prestable}
A nodal curve of genus $g$ over a scheme $T/k$ is a proper, flat
morphism $\pi \colon C \to T$ of finite type such that
for every geometric point $\Spec\ \eta$ of $T$, the fiber over $\Spec\ \eta$
is a connected curve of genus $g$ with only nodes as singularities.

Given $(g, S, \AA, \beta)$ as above, 
a prestable map of type $(g, \AA, \beta)$ over $T$ is a tuple
$(C, \pi, s, f)$ where
$\pi \colon C \to T$ is a nodal curve of genus $g$,
$s = (s_i)_{i \in S}$ is an $S$-tuple
of sections $s_i \colon T \to C$, and $f$ is a map
$f \colon C \to V$ with $f_*([C]) = \beta$, such that
\begin{enumerate}
\item
the image of any section $s_i$ with positive weight $\AA(i) > 0$
lies in the smooth locus of $C/T$,
\item
for any subset $I \subset S$ such that the intersection
$\bigcap_{i \in I} s_i(T)$ of the corresponding sections is non-empty, 
we have $\sum_{i \in I} \AA(i) \le 1$.
\end{enumerate}
\end{Def}

\begin{Def}					\label{def:stable}
A stable map of type $(g, \AA, \beta)$ over $T$ is a prestable map 
$(C, \pi, s, f)$ of the same type such that
$K_{\pi} + \sum_{i \in S} \AA(i) s_i + 3 f^*(M)$
is $\pi$-relatively ample for some ample divisor $M$ on $V$.
\end{Def}
We will often call such a curve $(g, \AA)$-stable when the homology class
$\beta$ is irrelevant.

\begin{Rem} \label{rem:geompoints}
Assume that $(C, \pi, s, f)$ is a $(g, \AA)$-prestable map over $T$.
Then it is $(g, \AA)$-stable if and only if it is $(g, \AA)$-stable
over geometric points of $T$.
\end{Rem}

Over an algebraically closed field, ampleness of 
$K_{\pi} + \sum_{i \in S} \AA(i) s_i + 3 f^*(M)$ can only fail on
irreducible components $C$ that are of genus 0
and get mapped to a point by $f$. Precisely, if
$n_C$ is the number of inverse images of nodal points in the normalization,
then ampleness is equivalent
to $n_{C} + \sum_{i \colon s_i \in C} \AA(i) > 2$.

In particular, stability can be checked with an arbitrary ample divisor
$M$; if all sections have weight 1 (we will write this as
$\AA = \underline{1}_S$), weighted stability agrees with the
definition of a stable map by Kontsevich.

We call the data $g, S, \AA, \beta$ {\it admissible}, if
$\beta \neq 0$ or $2g - 2 + \sum_{i\in S} \AA(i) > 0$, and if
$\beta$ is bounded by the characteristic
(cf. \cite[Theorem 3.14]{BehrendManin}: this means that $k$ has
characteristic zero, or that $\beta \cdot L < \Char\ k$ for some very
ample line bundle $L$ on $V$).
\begin{Thm}			\label{thm:DMstacks}
Given admissible data $g, S, \AA, \beta$,
let $\Mbar_{g, \AA}(V, \beta)$ be the category
of stable maps of type $(g, \AA, \beta)$ and their isomorphisms, with 
the standard structure as a groupoid over schemes over $\Spec\ k$.

This category is a proper algebraic Deligne-Mumford stack of
finite type.
\end{Thm}
The property of being a stack follows from standard arguments. The geometric
properties are proven in section \ref{subsect:geom-props}. Some of their proofs
are simplified by the use of the contraction morphism from the
Kontsevich moduli space $\Mbar_{g,S}(V, \beta)$ to the space of weighted stable
maps as discussed in the next section; hence their existence will be proved
first.

\subsection{Reduction morphisms for weight changes.}
If $\beta \neq 0$, consider the open and dense substack
$$ C_{g, S}(V, \beta) \subset \Mbar_{g, \AA}(V, \beta)
$$
of maps that do not
contract any irreducible component of genus zero, and for which all marked
sections are distinct.
By some abuse of language we will call $C_{g, S}(V, \beta)$
the ``configuration space''.
Since any such map is stable regardless of the choice of weights,
$C_{g, S}(V, \beta)$ does not depend on $\AA$.
Every $\Mbar_{g, \AA}(V, \beta)$ is a compactification of $C_{g, S}(V, \beta)$,
and thus all the moduli stacks for different $\AA$
are birational. The following proposition gives actual morphisms,
provided that the weights are comparable.
They will be analyzed in more detail in
\ref{sect:chamber-decomposition}.

Consider two weights
$\AA, \BB \colon S \to \Q \cap [0, 1]$
such that $\AA(i) \ge \BB(i)$ for all $i \in S$; we will just write
$\AA \ge \BB$ from now on. Any $(g, \AA)$-stable map is obviously
$(g, \BB)$-prestable, but it may not be $(g, \BB)$-stable. However, we can
stabilize the curve with respect to $\BB$:

\begin{Prop} 			\label{prop:reduction}
If $g, S, \beta, \AA \ge \BB$ are as above, there is a natural
reduction morphism
$$
\rho_{\BB, \AA} \colon \Mbar_{g, \AA}(V, \beta) \to \Mbar_{g, \BB}(V, \beta).
$$
It is surjective and birational.\footnote{It is an isomorphism over
the open subset $U := C_{g, S}(V, \beta)$,
which satisfies the following strong density property: for any
open subset $V$, there is no non-zero section $f \in \OO(V)$ that vanishes
on $U \cap V$; in other words, the complement is nowhere dense and does
not have \emph{additional} nilpotent structure.}
Over an algebraically closed field $\eta$, it is given by adjusting
the weights and then successively contracting all $(g, \BB)$-unstable
components. 

Given three weight data $\AA \ge \BB \ge \CC$, the reduction morphisms
respect composition:
$\rho_{\CC, \AA} = \rho_{\CC, \BB} \circ \rho_{\BB, \AA}$.
\end{Prop}

In particular, every moduli space $\Mbar_{g, \AA}(V, \beta)$ is a
birational contraction of the Kontsevich moduli space
$\Mbar_{g, S}(V, \beta) = \Mbar_{g, \underline{1}_S}(V, \beta)$.

\subsection{Proofs of the geometric properties.} \label{subsect:geom-props}

As in the case of $(g, \AA)$-stable curves, 
the following vanishing result is essential to ensure that all constructions
are compatible with base change:

\begin{Prop}\cite[Proposition 3.3]{Has03}  \label{prop:cohomology-vanishing}
Let $C$ be a connected nodal curve of genus $g$ over an algebraically closed
field, $D$ an effective divisor
supported in the smooth locus of $C$, and $L$ an invertible sheaf with
$L \cong \omega_C^k(D)$ for $k > 0$.

1. If $L$ is nef, and $L \neq \omega_C$, then $L$ has vanishing higher
cohomology.

2. If $L$ is nef and has positive degree, then $L^N$ is basepoint free
for $N \ge 2$.

3. If $L$ is ample, then $L^N$ is very ample when $N \ge 3$.

4. Assume $L$ is nef and has positive degree, and let $C'$ denote the image
of $C$ under $L^N$ with $N \ge 3$. Then $C'$ is a nodal curve with the same
arithmetic genus as $C$, obtained by collapsing the irreducible components
of $C$ on which $L$ has degree zero. Components on which $L$ has positive
degree are mapped birationally onto their images.
\end{Prop}

\subsubsection{Stability and geometric points.}
We will first show how remark \ref{rem:geompoints} follows from this
proposition:
Consider the line bundle
$$
L = \omega_C^k(k\sum_{i \in S} \AA(i) s_i) \otimes f^*(\OO(M))^{3k},$$
where $k$ is such that all numbers $k \AA(i)$ are integral. Then by the
proposition and the base change theorems, 
formation of $P := \Proj (\pi_* (L^{N}))$
commutes with base change. By definition, $L$ is relatively ample iff
the induced morphism $p \colon C \to P$
is defined everywhere and an open immersion. By
\cite[expos\'e I, Th\'eor\`eme 5.1]{SGA1}, this is the case if and only
if $p$ is everywhere defined, radical, flat and unramified.
All these conditions can be checked on
geometric fibers (for flatness, this follows from \cite[IV, Th\'eor\`eme
11.3.10]{EGA}, for unramifiedness from the conormal sequence).

\subsubsection{Reduction morphisms.} 		\label{proof-reduction}
By Grothendieck's descent theory, $\Mbar_{g, \AA}(V, \beta)$
is a stack in the \'etale topology, i.~e. the Isom functors are
sheaves and any \'etale descent datum is effective.
We first show the existence of the natural reduction morphisms
$\rho_{\BB, \AA}$ as maps between these abstract stacks.
This will enable us to use the results
of \cite{BehrendManin} on $\Mbar_{g,S}(V, \beta)$
to shorten our proofs.

Using the vanishing result \ref{prop:cohomology-vanishing},
the proof of proposition \ref{prop:reduction}
is completely analogous to that of theorem 4.1 in \cite{Has03}: Let
$\BB_\lambda = \lambda \AA + (1 - \lambda) \BB$, and let 
$1 = \lambda_0 > \lambda_1 > \dots > \lambda_N = 0$ be a finite set such
that for all $\lambda \not \in \{\lambda_0, \dots, \lambda_N \}$,
the following condition holds:
\begin{itemize}
\item
{\it There is no subset $I \subset S$ such that
$\sum_{i \in I} \BB_\lambda(i) = 1$ and
$\sum_{i \in I} \BB_1(i) \neq 1$.} (*)
\end{itemize}

We will construct $\rho_{\BB, \AA}$ as the composition
$\rho_{\BB, \AA} = \rho_{\BB(\lambda_N), \BB(\lambda_{N-1})} \circ \dots
\circ \rho_{\BB(\lambda_1), \BB(\lambda_0)} $. This means we can assume
that the condition (*) holds for all $0 < \lambda < 1$.

Fix an ample divisor $M$ on $V$, and fix a natural number $k$ so that
$k\BB(i)$ is an integer for all $i$. Let $L$ be the invertible sheaf
$L := \omega^k_C(k\sum_{i \in S} \BB(i) s_i) \otimes f^*(M)^3$ for any
$(g, \AA)$-stable map $f \colon C \to V$ over $T$. Due to
condition (*), it is nef; also it has positive degree. Let $C'$ be the
image of $C$ under the map induced by $L^N$ for some $N \ge 3$, i.e.
$C' = \Proj\ \RR$ where $\RR$ is the graded sheaf of rings on $T$ given
by $\RR_l = \pi_* ((L^N)^l)$. Let $t \colon C \to C'$ be the natural
map, and let $s_i' = t \circ s_i$.
By the same arguments as in the non-weighted case, $C'$ is a nodal
curve of genus $g$, and $s_i'$ lie in the smooth locus whenever $\BB(i) > 0$.
By proposition \ref{prop:cohomology-vanishing},
$L$ has vanishing higher cohomology; so 
the formation of $\pi_* ((L^N)^l)$ and hence that of $C'$ commutes with
base change.
Over an algebraically closed field, this morphism
agrees with the description via contraction of unstable components. In
particular, $C'$ is $(b, \BB)$-prestable.

The original $f$ factors via the induced morphism $f'\colon C' \to V$.
Let $L'$ be the line bundle
$L' := \omega^k_{C'}(k \sum_{i \in S} \BB(i) s_i) \otimes {f'}^*(M)^3$.
Then $t_* L = L'$; hence 
$L'$ is ample and  $(C', \pi', s', f')$ is a $(g, \BB)$-stable map.
The induced
morphism $T \to \Mbar_{g, \BB}(V, \beta)$ commutes with base change
and thus yields the map
$\rho_{\BB, \AA}$ between stacks as claimed.

To prove surjectivity, it is sufficient
to show that every $(g, \BB)$-stable map $(C, s, f)$ over
an algebraically closed field $K$ is 
the image of  some $(g, \AA)$-stable map $(C', s', f')$ over $K$.
It is obvious how to construct $C'$: 
If $I \subset S$ is a subset of the labels such that
condition (2) of definition \ref{def:prestable} is violated for
the weight data $\AA$, i.e.
the marked points $s_i, i \in I$ coincide and
$\sum_{i \in I} \AA(i) > 1$, we can attach a copy of $\P^1(K)$ at
this point, move the marked points to arbitrary but different points
on $\P^1$, and extend the map constantly along $\P^1$.

Birationality (for $\beta \neq 0$) follows from the fact that
$\rho_{\BB, \AA}$ is an isomorphism over the configuration space
$C_{g, S}(V, \beta)$, which is a dense and open subset.
The compatibility with composition follows immediately
once we have shown the the moduli spaces are separate: the two
morphisms $\rho_{\CC, \AA}$ and $\rho_{\CC, \BB} \circ \rho_{\BB, \AA}$
agree on the configuration space.

\begin{Prop}
The diagonal
$\Delta \colon \Mbar_{g, \AA}(V, \beta) \to
	\Mbar_{g, \AA}(V, \beta) \times
	\Mbar_{g, \AA}(V, \beta)$
is representable, separated and finite.
\end{Prop}

Let $(C_1, \pi_1, s_1, f_1)$ and
$(C_2, \pi_2, s_2, f_2)$ 
be two families of $(g, \AA)$-stable maps to $V$ over a scheme $T$. 
We have to show that
$\Isom((C_1, \pi_1, s_1, f_1),
(C_2, \pi_2, s_2, f_2))$ is represented by a 
scheme finite and separated over $T$. Since $V$ is projective and $\beta$ is
bounded by the characteristic,
we can use exactly the same argument as in the
proof of \cite[Lemma 4.2]{BehrendManin}: one shows that \'etale locally on $T$, one can
extend the set of labels to $S \cup S'$ and 
find additional $S'$-tuples of sections $(s_1)'$ and
$(s_2)'$, such that
$(C_1, \pi_1, s_1 \cup s'_1)$ and
$(C_2, \pi_2, s_2 \cup s'_2)$ are
$(g, \AA \cup \underline{1}_{S'})$-stable curves, and that there is
a natural closed immersion
$$
\Isom ((C_1, \pi_1, s_1, f_1), (C_2, \pi_2, s_2, f_2)) \to
\Isom ((C_1, \pi_1, (s_1, s_1')), (C_2, \pi_2, (s_2, s'_2))).
$$
Sine $\Mbar_{g, \AA \cup \underline{1}_{S'}}$ has a representable,
separated and finite diagonal by \cite{Has03}, the claim of the proposition
follows.

\subsubsection{Existence as Deligne-Mumford stacks.}
In particular, the diagonal is proper and thus the moduli stack separated.
As $\Mbar_{g, \underline{1}_S}(V, \beta)$ is proper
and the reduction morphism
$\rho_{\AA, \underline{1}_S} \colon \Mbar_{g, S}(V, \beta)
\to \Mbar_{g, \AA}(V, \beta)$ is surjective, 
$\Mbar_{g, \AA}(V, \beta)$ is also proper.

Finally, the existence of a flat covering of
finite type follows with almost the same argument as the one in
\cite{BehrendManin}, following Proposition 4.7 there. However,
some changes are required, so we spell it out in detail:
We write $\AA_n = \AA \cup \underline{1}_{\{1, \dots, n\}}$
for the weight data obtained from $\AA$
by adding $n$ weights of 1.
Let $\Mbar^o_{g, \AA_n}(V, \beta)$ be the open substack
of $\Mbar_{g, \AA_n}(V, \beta)$ where the
additional sections of weight one lie in the smooth locus of
$C_{g, \AA}(V, \beta)$ and away from the existing sections (in
other words, the open substack where the map is already $(g, \AA)$-stable 
after forgetting the additional sections).
The obvious forgetful map
$$
\phi_{\AA, \AA_n}^0 \colon \Mbar^o_{g, \AA_n}(V, \beta)
\to \Mbar_{g, \AA}(V, \beta)
$$
is smooth and in particular flat.
Let $U^0_{g, \AA_n}(V, \beta)$ be the open substack
of $\Mbar^o_{g, \AA_n}(V, \beta)$ where the curve
is already $(g, \AA_n)$-stable as a curve. Then for high enough $n$, the
restriction $\phi_{\AA, \AA_n}^o|_{U^0_{g, \AA_n}(V, \beta)}$
to this substack is surjective.
On the other hand, $U^0_{g, \AA_n}(V, \beta)$ is an open substack
of the relative morphism space
$\Mor_{\Mbar_{g, \AA_n}}(V, \beta)$ (parameterizing
maps $T \to \Mbar_{g, \AA_n}$ together with a map of the
pull-back of the universal curve $C_{g, \AA_n}$ to $V$).
So a flat presentation of the morphism space induces one for
$\Mbar_{g, \AA}(V, \beta)$.

\section{Elementary morphisms}
					\label{sect:elementary-morphisms}

\subsection{Gluing morphisms.}
As in the non-weighted case, we can glue curves at marked points, but 
to guarantee that the resulting curves are prestable, we have to assume that
both labels have weight 1:

Let $g_1, S_1, \AA_1, \beta_1$ and $g_2, S_2, \AA_2, \beta_2$ be
weight data, such that the
extensions $g_i, S_i \cup \{0\}, \AA_i \cup \{0 \mapsto 1\}, \beta_i$
by an additional label of weight 1 are admissible.
Denote by $\ev_0$ be the evaluation morphisms
$\ev_0 \colon \Mbar_{g_i, \AA_i \cup \{1\}}(V, \beta_i) \to V$ given
by evaluating the additional section: $\ev_0 = f \circ s_0$.
Similarly, let $g, S, \AA, \beta$ be weight data such that
$g, S \cup \{0, 1\}, \AA \cup \{1, 1\}, \beta$ is admissible, and
let $\ev_{0}, \ev_1$ be the additional evaluation morphisms.

\begin{Prop}					\label{prop:gluing}
There are natural gluing morphisms
\begin{align*}
\left(\Mbar_{g_1, \AA_1 \cup \{1\}}(V, \beta_1)
\times \Mbar_{g_2, \AA_2 \cup \{1\}}(V, \beta_2) \right)
\times_{V \times V} V
&\to \Mbar_{g_1 + g_2, \AA_1 \cup \AA_2}(V, \beta_1 + \beta_2) \\
\intertext{and}
\Mbar_{g, \AA \cup \{1, 1\}}(V, \beta) \times_{V \times V} V
& \to \Mbar_{g+1, \AA}(V, \beta).
\end{align*}
The product over $V \times V$ is taken via the morphism
$(\ev_0, \ev_0)$
respectively $(\ev_0, \ev_1)$ on the left, and the diagonal
$\Delta \colon V \to V \times V$ on the right.
\end{Prop}

There is nothing new to prove here, except to note that the weight of 1
guarantees that the marked sections (of positive weight) 
do not meet the additional node on the glued curve.

\begin{Prop-s}					\label{prop:target_change}
Let $\mu \colon V \to W$ be a morphism, and $(g, S, \AA, \beta)$ be admissible
data for $V$, such that $(g, S, \AA, \mu_*(\beta))$ is also admissible.
Then there is a natural push-forward
$$
\Mbar_{g, \AA}(V, \beta) \to \Mbar_{g, \AA}(W, \mu_*(\beta))
$$
that is obtained by composing the maps with $\mu$, followed by stabilization.
\end{Prop-s}

One could adapt the proof of \cite{BehrendManin} to the weighted case;
instead, we give a proof analogous to the one in section
\ref{proof-reduction}.

Let $f \colon C \to V$ be the universal map over
$\Mbar_{g, \AA}(V, \beta)$, let $f' = \mu \circ f$ be the induced
map to $W$, and let $M'$ be an ample divisor on $V'$.
By the assumptions, the divisor
$D' = K_\pi + \sum_{i \in S} \AA(i) s_i + 3 {f'}^*M'$ has positive degree;
however, it need not be nef. Hence we 
consider $D = K_\pi + \sum_{i \in S} \AA(i) s_i + 3 f^*M$ and
$D(\lambda) = \lambda D + (1 - \lambda) D'$ for $0 \le \lambda \le 1$.
Let $\{ \lambda_1, \dots, \lambda_N \}$ be the set of $\lambda$ for
which the degree of $D(\lambda)$ is zero on any irreducible component
of $C$, and let $k_r, r = 1\dots N$ be an integer such that $k_r \lambda_r$
and $k_r\AA(i), i \in S$ is integer.

Then
$L_1 = \omega^{k_1} (k_1 \sum_{i \in S} \AA(i) s_i + k_1 (3 f^*M \lambda_1 + (1 - \lambda_1) 3 {f'}^*M'))$
is a nef invertible sheaf on $C$ for which
proposition \ref{prop:cohomology-vanishing} applies. Hence $C_1$ defined by
$C_1 := \Proj\ \RR_1$ and $(\RR_1)_l = \pi_*(L_1^{3l})$ is again a flat
nodal curve of genus $g$, contracting all components
of $C$ on which $L_1$ fails to be ample, and $f'$ factors via a unique
morphism $f_1 \colon C_1 \to W$.
We proceed inductively to obtain $f_N \colon C_N \to W$ on which $D'$ is
ample; this induces the map of moduli stacks.
Note that $C \to C_N \to W$ is the universal factorization of $f'$ such
that $f_N \colon C_N \to W$ is a $(g, \AA)$-stable map.

\begin{Prop-s} 						\label{prop:forget}
Given admissible weight data $(g, S, \AA, \beta)$, let
$(g, S \cup \{*\}, \AA \cup \{a\} = \AA \coprod \{ * \mapsto a \}, \beta)$
 be the weight data
obtained by adding a label $\{*\}$ of arbitrary weight
$a \in \Q \cap [0, 1]$.
There is a natural forgetful map
$$
\phi_{\AA, \AA \cup \{a\}} \colon
\Mbar_{g, \AA \cup \{a\}}(V, \beta) \to \Mbar_{g, \AA}(V, \beta)
$$
obtained by forgetting the additional section and stabilization. If
$a = 0$, then $\phi_{A, A \cup \{0\}}$ is the universal curve over
$\Mbar_{g, \AA}(V, \beta)$.
\end{Prop-s}

We can construct this map as the composition
$$
\phi_{\AA, \AA \cup \{0\}} \circ {\rho_{\AA \cup \{0\}, \AA \cup \{a\}}} \colon
\Mbar_{g, \AA \cup \{a\}}(V, \beta)
\to \Mbar_{g, \AA \cup \{0\}}(V, \beta)
\to \Mbar_{g, \AA}(V, \beta).
$$
The second morphism $\phi_{\AA, \AA \cup \{0\}}$
is the naive forgetful morphism, as a map is
$(g, \AA \cup \{0\})$-stable if and only if it is
$(g, \AA)$-stable. 

\begin{Prop-s}					\label{prop:combine}
Let $S' \coprod S'' = S$ be a partition of the set of labels such that
$\AA(S'') = \sum_{i \in S''} \AA(i) \le 1$. Then there is a natural map
$$
\Mbar_{g, \AA|_{S'} \cup \{\AA(S'')\}}(V, \beta)
\to \Mbar_{g, \AA}(V, \beta).
$$
\end{Prop-s}
It is given by setting $s_i = s_*$ for all $i \in S'$. It identifies
$\Mbar_{g, \AA|_{S'} \cup \{* \mapsto \AA(S'')\}}(V, \beta)$ with the
locus of $\Mbar_{g, \AA}(V, \beta)$ where all $s_i, i \in S''$
agree.

\subsection{Weighted marked graphs.}	\label{sect:first-graphs-sketch}
A graph was defined in \cite{BehrendManin} as a quadruple
$\tau = (V_\tau, F_\tau, \partial_\tau, j_\tau)$ of a set of vertices
$V_\tau$, a set of flags $F_\tau$, a morphism
$\partial_\tau \colon F_\tau \to V_\tau$ and an involution
$j_\tau \colon F_\tau \to F_\tau$.
We think of a graph in terms of its geometric realization: it is obtained
by identifying in the disjoint union
$\coprod_{f \in F_\tau} [0,1]$ the points 0 for all
flags $f$ attached to the same vertex via $v = \partial_\tau(f)$, and
the points 1 for all orbits of $j_\tau$.
A flag $f$ with $j_\tau(f) = f$ is called a \emph{tail} of the vertex
$\partial_\tau(f)$, whereas a pair
$\{f, j_\tau(f)\}$ for $f \neq j_\tau(f)$ is called an \emph{edge},
connecting the (not necessarily distinct) vertices $\partial_\tau(f)$ and
$\partial_\tau(j_\tau(f))$. 

Given a projective variety $V$, a weighted modular $V$-graph
is a graph $\tau$ together with a genus $g \colon V_\tau \to \Z_{\ge 0}$,
a weight data $\AA \colon F_\tau \to \Q \cap [0,1]$ such that
$\AA(f) = 1$ for all flags that are part of an edge, and a marking
$\beta \colon V_\tau \to H_2^+(V)$. To any weighted stable map we can
associate its dual graph: a vertex for every irreducible component, an
edge for every node, and a tail for every marked section. Conversely, to every
weighted modular graph we can associate the moduli space of tuples
of weighted stable maps $f_v \colon C_v \to V$ of type
$(g(v), S_v = \{f \in F_\tau \colon \partial(f) = v \}, \AA|_{S_v}, \beta(v))$,
such that for every edge $\{f, f' = j_\tau(f)\}$ connecting the vertices
$v = \partial_\tau(f)$ and $v' = \partial_\tau(f)$, the corresponding
evaluation morphisms are identical: $f_v \circ s_f = f_{v'} \circ s_{f'}$.
Via gluing, this gives a single weighted stable map $f \colon C \to V$; if
all $C_v$ are smooth, its dual graph will give back $\tau$.

The moduli space $\Mbar_{g, \AA}(V, \beta)$ corresponds to the one-vertex
graphs with the set $S$ of tails.
The morphisms constructed in this section correspond to elementary morphisms
between graphs with one and two vertices.
Extending this set of morphisms to higher codimension boundary strata,
indexed by graphs with more vertices, naturally leads to a category of
weighted stable marked graphs.
We will adopt this viewpoint in
\ref{sect:graph-language}, and show that we get a functor
$\Mbar$ from the graph category to Deligne-Mumford stacks over $k$.

\section{Birational behavior under weight changes}
			\label{sect:chamber-decomposition}

For this section, we will fix $g, S, V, \beta$, and analyze more systematically
the reduction morphisms $\rho_{\AA, \BB}$ of proposition \ref{prop:reduction}
for varying weight data $\AA, \BB$. Assume that $g, V, \beta$ are such that
$\Mbar_{g, \AA}(V, \beta)$ is not empty.

\subsection{Exceptional locus and reduction morphism as blow-up}

\begin{Prop}\cite[Proposition 4.5]{Has03} 		\label{prop:boundary}
Assume we have weight data $\AA \ge \BB > 0$.
The reduction morphism $\rho_{\BB, \AA}$ contracts the boundary
divisors $D_{I, J}$ given as the image of the gluing morphism
$$
\Mbar_{0, \AA|_I \cup \{1\}}(V, 0) \times_V
\Mbar_{g, \AA|_J \cup \{1\}}(V, \beta)
\to \Mbar_{g, \AA}(V, \beta)
$$
for all partitions $I \coprod J = S$ of $S$
with
$$
\sum_{i \in I} \AA(i) > 1 \quad \text{and} \quad
b_I := \sum_{i \in I} \BB(i) \le 1.
$$
There is a factorization of
$\rho_{\BB, \AA}|D_{I, J}$ via
$$
\Mbar_{0, \AA|_I \cup 1}(V, 0) \times_V
\Mbar_{g, \AA|_J \cup 1}(V, \beta)
\to \Mbar_{g, \AA|_J \cup \{1\}}(V, \beta)
\to \Mbar_{g, \AA|_J \cup \{b_I\}}(V, \beta).
$$
\end{Prop}

We may assume that there is just one such $I$
and that $b_I = 1$. The stabilization contracts
components on which
$\omega^k_C(k\sum_{i \in S} \BB(i) s_i) \otimes f^*(M)^3$
has degree zero. Such a component
can only be a smooth irreducible component of genus zero
that is mapped to a point, meets the other components in a single node
and contains exactly those marked sections $s_i$ with $i \in I$.

In particular, the exceptional locus of $\rho_{\BB, \AA}$ is given
by all $D_{I, J}$ for partitions $I \cap J = S$ as above with the additional
condition $\abs{I} > 2$.
When all sets $I \subset S$ such that
$\sum_{i \in I} \AA(i) > 1$ and $\sum_{i \in I} \BB(i) \le 1$ satisfy
$\abs{I} = 2$, then $\rho_{\BB, \AA}$ is an isomorphism.

\begin{Rem}					\label{rem:reduction-is-blowup}
Assume that for $\AA > \BB > 0$,
there is exactly one partition $I \coprod J = S$ of $S$ as
in the proposition. Then $\rho_{\BB, \AA}$ is the blow-up of
$\Mbar_{g, \BB}(V, \beta)$ along the substack
$C_{IJ} \cong \Mbar_{g, \BB|_J \cup \{b_I\}}(V, \beta)$ of weighted stable
curves where all section $s_i$ for $i \in I$ are identical.
\end{Rem}

We first show that there is a natural map from $\Mbar_{g, \AA}(V, \beta)$ to
the blow-up:
The divisor $D_{I, J}$ is the scheme-theoretic inverse image of $C_{IJ}$.
Further, it is a Cartier divisor: if $C$ is the universal curve over
$\Mbar_{g, \AA}(V, \beta)$, and $C'$ the pull-back of the universal curve
over $\Mbar_{g, \BB}(V, \beta)$, then $D_{I, J}$ is the zero locus of
the natural map $s_{i_0}^* \Omega_C \to s_{i_0}^* \Omega_{C'}$ of
the pull-backs of the relative cotangent sheaves for some
$i_0 \in I$.
By the universal property of blow-ups, this shows that
$\rho_{\BB, \AA}$ factors via the blow-up
$\rho' \colon M \to \Mbar_{g, \BB}(V, \beta)$ of $\Mbar_{g, \BB}(V, \beta)$
at $C_{IJ}$.

We now construct the inverse map. Let $C'$ be the pull-back of the universal
curve
along $\rho'$, let $E$ be the exceptional divisor of $\rho'$, and write
$\rho'^{-1}s_i \colon M \to C'$ for the pull-back of the sections
$s_i$ over $\Mbar_{g, \BB}(V, \beta)$.
Let $C_0$ be the common image $(\rho'^{-1}s_i)(E)$ of the exceptional
divisor for any $i \in I$, and let $C$ be blow-up of $C'$ at
$C_0$.
The center $C_0 \subset C'$ is a codimension two regular embedding, and
embeds as a Cartier divisor in
both $(\rho'^{-1}s_i)(M)$ for any $i \in I$, and in the restriction of $C'$
to $E$. Thus the fibers of $C$ over $E$ are obtained from that of the
universal curve over $\Mbar_{g, \BB}(V, \beta)$ by attaching a projective
line at the marked point given by any $s_i$ for $i \in I$, and
every section $\rho'^{-1}s_i$ lifts to a section $s_i \colon M \to C$
via the proper transform of $(\rho'^{-1}s_i)(M)$.

Over $E$, the image is contained in the attached projective line, away from
the node, as $s_i(M)$ and the fiber over $E$ meet transversely in $C'$.
Also, since the images of $s_i, i \in I$ intersect transversely in the
universal curve over $\Mbar_{g, \BB}(V, \beta)$, any tangent vector at a
point of $C_0$ tangent to all the images of $(\rho'^{-1}s_i)(M), i \in I$
is already
tangent to $C_0$; thus the sections $s_i \colon M \to C$ cannot all be
mapped to the same point of the projective line.

Hence, with the induced map to $V$, we have constructed a $(g, \AA)$-stable
map, and so a map $M \to \Mbar_{g, \AA}(V, \beta)$; it is an inverse
to the map in the opposite direction constructed above, as this is true
over $C_{g, S}(V, \beta)$ and both stacks are separated.

\begin{Prop}
Let $\AA, \BB$ as in proposition \ref{prop:boundary}, except we allow
some weights of $\BB$ to be zero.
Let $i \in S$ be a label with $\AA(i) > \BB(i) = 0$.
Then $\rho_{\AA, \BB}$
additionally contracts the boundary components
$C_{(g_1, 0, g_2), (I_1, I_0, I_2), (\beta_1, 0, \beta_2)}$ which are 
defined as the image of the gluing morphisms
$$
\Mbar_{g_1, \AA|_{I_1} \cup \{1\}}(V, \beta_1)
\times_V \Mbar_{0, \AA_{I_0 \cup \{i\} } \cup \{1, 1\}}(V, 0)
\times_V \Mbar_{g_2, \AA|_{I_2} \cup \{1\}}(V, \beta_2)
$$
$$
\to \Mbar_{g, \AA}(V, \beta)
$$
for all $g_1 + g_2 = g$, $\beta_1 + \beta_2 = \beta$ and disjoint
partitions $I_1 \cup I_0 \cup \{i\} \cup I_2 = S$ such that 
$\AA(j) = 0$ for $j \in I_0$.

The restriction $\rho_{\BB, \AA}$ factors via the projection of the
second component to a point.
\end{Prop}
In other words, this is the boundary component
of singular curves such that the section $s_i$ is contained in a node after
stabilization.

\subsection{Chamber decomposition}
We now assume $\beta \neq 0$, and consider the set of positive
weights $\DD_{n} = (0, 1]^S \subset \R^S$.
The walls $\WW_c$ and $\WW_f$ of the coarse and fine chamber decomposition,
respectively, are given by:\footnote{
The conditions $|S| < n-2$ and $|S| \le n-2$ for the
coarse and fine chamber decompositions, respectively, in \cite[section 5]{Has03}
are correct only when $g = 0$ and don't apply in our case as
we assumed $\beta \neq 0$.}
$$
\WW_c = \stv{\sum\nolimits_{i \in I} \AA(i) = 1} {I \subset S, 2 < |I|}
$$
$$
\WW_f = \stv{\sum\nolimits_{i \in I} \AA(i) = 1} {I \subset S, 2 \le |I|}
$$
Coarse and fine chambers are connected component of the complements
$\DD_n \setminus \WW_c$ and $\DD_n \setminus \WW_f$, respectively.

\begin{Prop} (cf. \cite[Proposition 5.1]{Has03})
The coarse chamber decomposition is the coarsest decomposition such
that $\Mbar_{g, \AA}(V, \beta)$ is constant in each chamber.
The fine chamber decomposition is the coarsest decomposition such
that the universal curve $\CC_{g, \AA}(V, \beta)$ is constant in each
chamber.
\end{Prop}

\begin{Cor}					\label{cor:universal-curve}
Let $\AA$ be positive weight data in the interior of a fine open chamber.
Then for small $\epsilon > 0$, the forgetful morphism
$\phi_{\AA, \AA \cup \{\epsilon\}}$ identifies
$\Mbar_{g, \AA \cup \{\epsilon\}}(V, \beta)$ with the universal curve
$C_{g, \AA}(V, \beta) \to \Mbar_{g, \AA}(V, \beta)$.
\end{Cor}
This holds by definition for $\epsilon = 0$, and 
it follows easily from the earlier propositions that
$\rho_{\AA \cup \{0\}, \AA \cup \{\epsilon\}}$ is an isomorphism.

\section{Virtual fundamental classes and Gromov-Witten invariants}
						\label{sect:virtual}

From now on, we assume additionally that the target $V$ is smooth.

\subsection{Expected properties}
\label{sect:expectedproperties}

The crucial step in the construction of Gromov-Witten invariants
is the construction of virtual fundamental classes of expected dimension:
$$
[\Mbar_{g, \AA}(V, \beta)]^{\virt} \in
A_{(1-g)(\dim V - 3) - K_V \cdot \beta + \abs{S}}
\Mbar_{g, \AA}(V, \beta)
$$
We will provide now a basic list of properties that such a
construction should satisfy, and proceed to draw some conclusions about
Gromov-Witten invariants in the remainder of the section. 

\begin{enumerate}
\item \emph{Mapping to a point.} 			\label{cl:point}
If $\beta = 0$, then 
$$
[\Mbar_{g, \AA}(V, 0)]^\virt = c_{g \dim V} (R^1 \pi_* f^* TV)
$$
\item \emph{Forgetting a tail.} 			\label{cl:forget}
Assume $\AA$ and $\epsilon$ are as in corollary \ref{cor:universal-curve},
so that $\phi_{\AA, \AA \cup \epsilon}$ is the universal curve over
$\Mbar_{g, \AA}(V, \beta)$.
In particular, this implies that $\phi_{\AA, \AA \cup \{\epsilon\}}$ is
flat, and thus defines a pull-back in intersection theory. We require
$$
\phi_{\AA, \AA \cup \epsilon}(V, \beta)^*
[\Mbar_{g, \AA}(V, \beta)]^\virt = 
[\Mbar_{g, \AA \cup \epsilon}(V, \beta)]^\virt.
$$

\item \emph{Combining tails.}				\label{cl:combine}
Assume we are in the situation of proposition \ref{prop:combine}.
Since all sections lie in the smooth locus of the curve,
$\mu_{S/S'}$ is a regular embedding, and we require that
$$
\mu_{S/S'}^! [\Mbar_{g, \AA}(V, \beta)]^\virt
= [\Mbar_{g, \AA|_{S'} \cup \{\AA(S'')\}}(V, \beta)]^\virt.
$$

\item \emph{Gluing.}				\label{cl:gluing}
We fix $g_1, S_1, \AA_1, g_2, S_2, \AA_2$ and some $\beta \in H_2^+(V)$.
Set $g = g_1 + g_2$ and $\AA = \AA_1 \cup \AA_2$. Consider the gluing
morphisms
$$
\mu_{\beta_1, \beta_2} \colon
\Mbar_{g_1, \AA_1 \cup \{1\}}(V, \beta_1) \times
\Mbar_{g_2, \AA_2 \cup \{1\}}(V, \beta_2) \times_{V \times V} V
$$
$$
\to
\Mbar_{g, \AA}(V, \beta)
$$
of proposition \ref{prop:gluing} for all $\beta_1, \beta_2$ with
$\beta_1 + \beta_2 = \beta$. The union of their images is the boundary
component in $\Mbar_{g, \AA}(V, \beta)$ given as the pull-back
$$
\xymatrix{
\Mbar_{(g_1, \AA_1) | (g_2, \AA_2)}(V, \beta) \ar[r] \ar[d] &
\Mbar_{g, \AA}(V, \beta) \ar[d] \\
\Mbar_{g_1, \AA_1 \cup \{1\}} \times \Mbar_{g_2, \AA_2 \cup \{1\}}
\ar[r]^{\hphantom{movemu}\mu} &
\Mbar_{g, \AA}
}
$$
Since the moduli spaces of weighted stable curves are smooth, 
$\mu$ is a l.c.i. morphism and defines a pull-back
$\mu^!  [\Mbar_{g, \AA}(V, \beta)]^\virt$.
On the other hand, via the diagonal $\Delta \colon V \to V \times V$, we
can pull-back the virtual fundamental class on the product
$\Mbar_{g_1, \AA_1 \cup \{1\}}(V, \beta_1)
	        \times \Mbar_{g_2, \AA_2 \cup \{1\}}(V, \beta_2)$
to the fiber product that is the source of $\mu_{\beta_1, \beta_2}$.
We require
$$
\sum_{\beta_1 + \beta_2 = \beta}
{\mu_{\beta_1, \beta_1}}_*
\Delta^! \left([\Mbar_{g_1, \AA_1 \cup \{1\}}(V, \beta_1)]^\virt
                \times [\Mbar_{g_2, \AA_2 \cup \{1\}}(V, \beta_2)]^\virt\right)
$$ $$
= \mu^!  [\Mbar_{g, \AA}(V, \beta)]^\virt.
$$

\item \emph{Kontsevich-stable maps.}			\label{cl:konts}
If all weights are 1, then
$[\Mbar_{g, \AA}(V, \beta)]^\virt$ agrees with the definition of
virtual fundamental classes of \cite{BF, Behrend:GW-alg}.

\item \emph{Reducing weights.}				\label{cl:reduct}
Given two set of weights $\AA > \BB$, we require compatibility
with the reduction morphism $\rho_{\BB, \AA}$:
$$
{\rho_{\BB, \AA}}_* [\Mbar_{g, \AA}(V, \beta)]^\virt
= [\Mbar_{g, \AA}(V, \beta)]^\virt
$$
\end{enumerate}

Evidently, properties (\ref{cl:point}), (\ref{cl:forget}) and
(\ref{cl:gluing}) are direct generalizations of properties satisfied by
the virtual fundamental classes of the non-weighted moduli spaces, while
(\ref{cl:combine}) and (\ref{cl:reduct}) are new.

\begin{Thm}						\label{thm:virt_class}
There is a system of virtual fundamental classes satisfying all of the
above properties.
\end{Thm}

While the Behrend-Fantechi construction can be applied to our situation
and provides virtual fundamental classes, we instead use (\ref{cl:konts})
and (\ref{cl:reduct}) as a definition, and prove that these classes
automatically satisfy the desired properties.

We postpone the proof of the above properties to
\ref{sect:virtual-and-graphs}, after having generalized them to a bigger
class of morphisms labelled by a category of weighted stable graphs. In the
remainder of the section we will instead proceed to give some consequences of
theorem \ref{thm:virt_class}.

\subsection{Gromov-Witten invariants}

As in the non-weighted case, one defines the $n$-point Gromov-Witten
invariants of $V$ depending on weights
$\AA \colon \{1, \dots, n\} \to [0, 1] \cap \Q$ via
\begin{align*}
\langle \hphantom{\gamma} \rangle_{g, \AA, \beta}
 \colon 
H^*(V)^{\otimes n} & \to \C
\\
\langle \gamma_1 \otimes \dots \otimes \gamma_n \rangle_{g, \AA, \beta}
& =  \int_{[\Mbar_{g, \AA}(V, \beta)]^\virt}
	\ev_1^*(\gamma_1) \cup \dots \cup \ev_n^*(\gamma_n)
\\
\intertext{
and Gromov-Witten invariants including gravitational descendants via}
\langle \tau_1^{k_1} \gamma_1 \cdots \tau_n^{k_n} \gamma_n
\rangle_{g, \AA, \beta}
& =  \int_{[\Mbar_{g, \AA}(V, \beta)]^\virt}
	\psi_1^{k_1} \ev_1^*(\gamma_1)
	\cup \dots \cup \psi_n^{k_n} \ev_n^*(\gamma_n)
\end{align*}
where $\psi_i$ is the tautological class associated to the section
$s_i$: $\psi_i = c_1(s_i^* \Omega_C)$ where $\Omega_C$ is the
relative cotangent bundle of the universal curve $C$ over
$\Mbar_{g, \AA}(V, \beta)$.

\begin{Prop} 		\label{prop:GW-without-gravity}
Gromov-Witten invariants without gravitational descendants do not depend
on the choice of weights $\AA$.
\end{Prop}
It is enough to prove this for two weights $\AA > \BB$. The evaluation
morphisms $\ev_i \colon \Mbar_{g, \AA}(V, \beta) \to V$ factor via the
reduction morphism $\rho_{\BB, \AA}$. Hence the claim follows
from property (\ref{cl:reduct}) and the projection formula.

\subsection{Extended modular operad}

Let $\AA_{m, n}$ be the weight data consisting of $m$ weights of one,
and $n$ weights of $\epsilon \le \frac 1n$.
The moduli spaces $\Mbar_{g, \AA_{m, n}}$ were called
$L_{g, m, n}$ in \cite{LosevManin-Ext} and studied more closely in
\cite{Manin-LgS}. Markings with weight one and $\epsilon$ are
white and black points in the language of \cite{LosevManin-Ext},
respectively: white points may not coincide with any other point, whereas
any number of black points are allowed to coincide.
Similarly, we write $L_{g, m, n}(V, \beta)$ for the moduli spaces of
weighted stable maps
$L_{g, m, n}(V, \beta) = \Mbar_{g, \AA_{m, n}}(V, \beta)$.

In  \cite{LosevManin-Ext}, the notion of an $\LL$-algebra was
introduced by a combinatorial description. It is an extension of the
graph-level description of the genus zero-part of a cohomological
field theory in the sense of \cite{KoMa}. Equivalently, it can be given
by a system of 
cohomology classes in $L_{0, m, n}$ (rather than classes in
$L_{0, m, 0} \cong \Mbar_{0, m}$). An $\LL$-algebra yields a (formal)
solution to the commutativity equations, which are extension of the
WDVV equations.

By the results of \cite{Manin-LgS}, the ''economy class description`` of an
$\LL$-algebra give in \cite[section 4.2.2]{LosevManin-Ext} can be translated
into the following cohomological description:

Let $(T; F, (,))$ be a triple consisting of two $\Z_2$-graded $\Q$-vector
spaces $T, F$, where the latter is equipped with a (super)symmetric
non-degenerate scalar product $(,)$. An $\LL$-algebra on
$(T; F, (,))$ over a $\Q$-algebra $R$ can be given as a collection of maps
\[
I_{0; m, n} \colon
T^{\otimes n} \otimes F^{\otimes m}
\to H_*(L_{0; m, n}) \otimes_\Q R
\]
being compatible with gluing of black points and the trace on $F$.

We obtain the \emph{$\LL$-algebra of quantum cohomology
of $V$ including gravitational descendants} as follows: Let $F = H^*(V,
\Q)$, equipped with the Poincar\'e pairing, and let $T = \bigoplus_{k \ge
0} z^k F$.  We denote by $\ev_1^W, \dots, \ev_m^W$ and $\ev_1^B, \dots,
\ev_n^B$ the evaluation maps $L_{0; m, n}(V, \beta) \to V$
induced by the marked sections of weight one and
$\epsilon$, respectively, and by
$\pi \colon L_{0; m, n}(V, \beta) \to L_{0; m, n}$ the forgetful map.
Let $\psi_i, i = 1 \dots n$ be the tautological
classes associated to the section $s_i^B$ of weight $\epsilon$.
Let $\Q[[q]]$ be the Novikov
ring of $V$, i.e. the formal completion of the polynomial ring over the
semigroup of effective classes in $H_2(V)/\text{torsion}$.

Then we define $I_{0; m, n}$ as
$$
I_{0; m, n} \left(
z^{k_1}\gamma_1 \otimes \dots \otimes z^{k_n} \gamma_n
\otimes \delta_1 \otimes \dots \otimes \delta_m
\right)
$$
$$
= \sum_{\beta \in H_2^+(V)}
q^\beta \PP\left(
\pi_* \left(
\prod_{i = 1}^n (\ev_i^B)^* \gamma_i \psi_i^{k_i} 
\prod_{i = 1}^n (\ev_i^W)^* \delta_i
	\cap [L_{0, m, n}(V, \beta)]^\virt \right) \right) 
$$
where $\pi \colon L_{0; m, n}(V, \beta) \to L_{0; m, n}$ is the forgetful
map, and $\PP(s) \in H^*L_{0; m, n}$ is the Poincar\'e dual of
$s \in H_*L_{0; m, n}$.

\begin{Thm}
The above definition of $I_{0; m, n}$ yields a cyclic $\LL$-algebra (in the
sense of the economy class description in \cite[section
4.2.2]{LosevManin-Ext}).
\end{Thm}

The only thing to check is the compatibility with gluing, in the formal
sense of \cite[diagram (4.8)]{LosevManin-Ext}. 
This holds due to property (4) of section \ref{sect:expectedproperties}.

\subsection{Comments.}
In \cite{LosevManin-Ext}, it was shown that the datum of an $\LL$-algebra
is equivalent to a geometric structure, a solution of the so-called
commutativity equation.
However, the structure of an $\LL$-algebra does not capture the complete
structure we have available:
\begin{enumerate}
\item
By property (6), the inclusion $F = z^0 F \subset T$ is compatible with
the reduction morphisms $L_{0, m, n} \to L_{0, m-1, n+1}$ in the obvious
sense.
\item 
Relating the gravitational descendants to intersection numbers
in $L_{0; m, n}$ by an analysis analogous to the one in \cite{KM-gravitation}
will, of course, lead to many more relations among the correlators.
\end{enumerate}
One might hope that these can be integrated in the geometric picture of
\cite{LosevManin-Ext}.

As a side remark, it is worth pointing out that the tautological
classes $\psi_i, i = 1\dots n$ in $L_{0; m, n}(V, \beta)$ are compatible with
pull-back along the forgetful morphism $L_{0; m, n+1}(V, \beta)$; this is
not true in the non-weighted case. 

\section{Graph-language}
						\label{sect:graph-language}

\subsection{Weighted marked graphs.}
						\label{section:weighted-graphs}
The elementary morphisms described in
\ref{sect:elementary-morphisms} 
generate a larger system of morphisms.
They are best modelled over a category of weighted marked graphs;
this category generalizes the category of marked graphs introduced in
\cite{BehrendManin} by introducing weights of tails. We follow
\cite[section 1]{BehrendManin} closely.

We recall from section \ref{sect:first-graphs-sketch} the definition of
a graph:
\begin{Def}\cite[Definition 1.1]{BehrendManin}
A graph $\tau$ is a quadruple $(F_\tau, V_\tau, j_\tau, \partial_\tau)$
of a finite set $V_\tau$ of vertices, a finite set $F_\tau$ of flags,
an involution $j_\tau \colon F_\tau \to F_\tau$ and a map
$\partial_\tau \colon F_\tau \to V_\tau$. We call
$S_\tau = \{ f \in F_\tau | j_\tau f = f \}$ the set of tails, and
$E_\tau = \{ \{f, j_\tau f \} | f \in F_\tau \ \text{and}\ j_\tau f \neq f \}$
the set of edges.
\end{Def}

\begin{Def}
A weighted modular graph is a graph
$\tau = (F_\tau, V_\tau, j_\tau, \partial_\tau)$
endowed with two maps
$g_\tau \colon V_\tau \to \Z_{\ge 0}$ and
$\AA_\tau \colon F_\tau \to \Q \cap (0, 1]$ such that
$\AA_\tau(f) = 1$ for all flags $f$ that are part of an edge, i.e. for which
$j_\tau(f) \neq f$.
\end{Def}
The number $g_\tau(v)$ is called the genus of a vertex, and
$\AA_\tau(f)$ the weight of a flag.

\begin{Def}			\label{def:combinat-m}
Given a semigroup $A$ with indecomposable zero, a
weighted $A$-graph $(\tau, \alpha)$
is a weighted modular graph $\tau$ with a map $\alpha \colon V_\tau \to A$.
A weighted marked graph is a pair $(A, (\tau, \alpha))$ where $A$
is a semigroup with indecomposable zero, and $(\tau, \alpha)$ is an $A$-graph.
\end{Def}
We will often omit $\alpha$ from the notation and call $\tau$ an $A$-graph.

Morphisms in the category of weighted marked graphs are generated by
two different types, \emph{combinatorial morphisms} and \emph{contractions}.
More precisely, since the associated geometric morphisms are contravariant
with respect to the combinatorial morphisms, and covariant with respect
to contractions, the morphisms will be generated by contractions and
formal inverses of the combinatorial morphisms.

Only condition (2) of the definition of a combinatorial morphism
of modular graphs (\cite[Definition 1.7]{BehrendManin}) needs to be adapted to
our situation:
\begin{Def}
Let $(\sigma, \alpha)$ and $(\tau, \beta)$ be weighted $A$-graphs.
A combinatorial
morphism $a \colon (\sigma, \alpha) \to (\tau, \beta)$ is a pair of maps
$a_F \colon F_\sigma \to F_\tau$ and $a_V \colon V_\sigma \to V_\tau$,
satisfying the following conditions:

\begin{enumerate}
\item The morphisms commute with $\partial$, i.e. we have
$a_V \circ \partial_\sigma = \partial_\tau \circ a_F$.
In particular, for any $v \in V_\sigma$ and
$w = a_V(v) \in V_\tau$,
we get an induced map $a_{V, v} \colon F_\sigma(v) \to F_\tau(w)$.

\item Consider the above map $a_{V, v}$. Then for any $f \in F_\tau(w)$,
the inequality
$$
\sum_{f' \in F_\sigma(v) \colon a_{V, v}(f') = f} \AA_\sigma(f') 
\le \AA_\tau(f)
$$
is satisfied.

\item
Let $\{f, \bar f\}$ be an edge of $\sigma$, i.e. $f \in F_\sigma, \bar f =
j_\sigma(f) \neq f$. Then there exist $n \ge 1$ and $n$ edges
$\{f_1, \bar f_1\}, \dots, \{f_n, \bar f_n\}$ of $\tau$ such that
$v_i := \partial_\tau(\bar f_i) = \partial_\tau(f_{i+1})$ and
$\beta(v_i) = 0$ for all $1 \le i < n$.

\item For every $v \in V_\sigma$ we have $\alpha(v) = \beta(a_V(v))$.

\item For every $v \in V_\sigma$ we have $g(v) = g(a_V(v))$.
\end{enumerate}

A {\it combinatorial morphism of weighted marked graphs} 
$(B, \sigma, \beta) \to (A, \tau, \alpha)$ is a pair $(\xi, a)$ where
$\xi \colon A \to B$ is a homomorphism of semigroups, and
$a \colon (\sigma, \beta) \to (\tau, \xi \circ \alpha)$ is
a combinatorial morphism of $B$-graphs.
\end{Def}

Note that we do not require that $j_\sigma$ and $j_\tau$ commute with 
$a_F$ and $a_V$; in particular, $\sigma$ could be obtained from $\tau$
by cutting an edge into two tails. Other examples of combinatorial morphisms
are morphisms adding tails or adding connected components.
There are essentially two new types of morphisms compared to the non-weighted
case:
\begin{enumerate}
\item
(\emph{Combining tails.}) Consider a subset $\{t_1, \dots, t_n\} \in
F_\sigma(v)$ of tails attached to a vertex $v$, and assume that its sum of
weights satisfies $\sum_i \AA_\sigma(t_i) \le 1$. Then we can form a new
graph $\tau$ by replacing the tails $\{t_1, \dots, t_n\}$ with a single
tail $\bar t$ of weight $\AA_\tau(\bar t) := \sum_i \AA_\sigma(t_i)$. 

\item
(\emph{Increasing the weights.}) This means that $(\tau, \beta)$ are identical
to $(\sigma, \alpha)$ as modular graphs, but the weight data $\AA_\tau$
satisfies $\AA_{\tau} \ge \AA_\sigma$.
\end{enumerate}

We refer to \cite[Definition 1.3]{BehrendManin} for 
the definition of a contraction $\phi \colon \tau \to \sigma$
of graphs. It is obtained by collapsing 
a subgraph consisting entirely of edges (and the adjoining vertices) to
one vertex for every connected component of the subgraph. 
It is given by an injective map
$\phi^F \colon F_\sigma \to F_\tau$ (which is bijective on tails) and
a surjective map $\phi_V \colon V_\tau \to V_\sigma$.

\begin{Def} 					\label{def:contraction}
A contraction of weighted marked graphs
$\phi \colon (\tau, \beta) \to (\sigma, \alpha)$
is a contraction $\phi \colon \tau \to \sigma$ of graphs such that
\begin{enumerate}
\item
$\alpha(v) = \sum_{w \in \phi_V^{-1}(v)} \alpha(w)$ for all $v \in V_\sigma$,
\item
$g(v) = \sum_{w \in \phi_V^{-1}(v)} \alpha(w) + H^1(\abs{\tau_v})$
for all $v \in V_\sigma$ and $\tau_v$ being the subgraph of $\tau$ being
collapsed onto $v$, and
\item
$\AA_\tau(\phi^F(f)) = \AA_\sigma(f)$ for all tails $f \in S_\sigma$.
\end{enumerate}
\end{Def}

\begin{Def}
A vertex $v$ of a weighted modular $A$-graph $(\tau, \alpha)$
is called stable if $\alpha(v) \neq 0$ or
$2g(v) - 2 + \sum_{f \in F_\tau \colon \partial_\tau(f) = v} \AA_\tau(f) > 0$.
A graph is stable if all its vertices are stable.
\end{Def}

\begin{Rem} 
Let $(\tau, \alpha)$ be a weighted $A$-graph. There is a
unique weighted stable $A$-graph $(\tau^s, \alpha^s)$ and a combinatorial
morphism $(\tau^s, \alpha^s) \to (\tau, \alpha)$,
such that every combinatorial morphism $(\sigma, \beta) \to (\tau, \alpha)$
from a stable $A$-graph $(\sigma, \beta)$ factors uniquely
through $(\tau^s, \alpha^s)$.
\end{Rem}

The graph $(\tau^s, \alpha^s)$ is called the \emph{stabilization} of
$(\tau, \alpha)$. Similarly, there is a stabilization of weighted
modular graphs. The stabilization $\tau^s$ of the underlying modular
graph $\tau$ of an $A$-graph $(\tau, \alpha)$ is also called the
\emph{absolute stabilization}.

The stabilization $(\tau^s, \alpha^s)$ can be constructed via a
sequence of steps as below, following
\cite[Proposition 1.13]{BehrendManin}:
\begin{enumerate}
\item
If there is a connected component of $\tau$ that has only one vertex, and
this vertex is unstable, we remove this connected component from $\tau$.
\item
If there is an unstable vertex $v$ attached to one edge
$\{f_0, \bar f_0 = j_\tau(f_0)\}$
with $\partial_\tau(f_0) = v$, $\partial_\tau(\bar f_0) \neq v$
and $n \ge 0$ tails $f_1, \dots, f_n$, we remove the vertex $v$ and the
flags $f_0, \dots, f_n$ from the graph and modify $j$ such that
$j(\bar f_0) = \bar f_0$, i.e. the edge becomes a tail at the vertex
$\partial_\tau(\bar f_0)$ with weight one.
\item
If there is an unstable vertex $v$ attached to two edges 
$\{f_1, \bar f_1 = j_\tau(f_1)\}$ and $\{f_2, \bar f_2 = j_\tau(f_2)\}$ with
$\partial_\tau(f_i) = v$ and $\partial_\tau(\bar f_i) \neq v$, we
remove $v$ and the tails $f_i$ from the graph, and modify $j$ such that
$j(\bar f_1) = \bar f_2$. In other words, we combine the tails
$\bar f_1, \bar f_2$ to form a new edge.
\end{enumerate}
At every step, any combinatorial morphism $(\sigma, \beta) \to
(\tau, \alpha)$, where $(\sigma, \beta)$ is a stable $V$-graph, factors
uniquely through the new graph, and the claim of the remark follows by
induction on the number of unstable vertices.

\begin{Def}
Let $(A, \tau)$ and $(B, \sigma)$ be weighted stable marked graphs. A
morphism $(A, \tau) \to (B, \sigma)$ is quadruple
$(\xi, a, \tau', \phi)$ where $\xi \colon A \to B$ is a homomorphism of
semigroups, $\tau'$ is a weighted stable $B$-graph, $a \colon \tau' \to \tau$
makes $(\xi, a)$ into a combinatorial morphism of weighted marked graphs, and
$\phi \colon \tau' \to \sigma$ is a contraction of $B$-graphs.
\end{Def}
$$
\xymatrix{
B & \tau' \ar[r]^\phi \ar[d]^a & \sigma \\
A \ar[u]^\xi & \tau
}
$$
We think of this morphism as the composition of
$\phi$ with the inverse of $(\xi, a)$, except that $(\xi, a)$ itself
is not a morphism in the category of weighted stable marked graphs. As
explained earlier, this construction is motivated by the fact that the
geometric morphisms are covariant with respect to contractions, but
contravariant with respect to combinatorial morphisms.

\label{stable-pullback}
To define compositions, we need the definition of stable pullback; the
construction of \cite{BehrendManin} applies with minor changes.
Given a combinatorial morphism of weighted marked graphs
$(a, \xi) \colon (B, \rho) \to (A, \tau)$ and a contraction of
weighted $A$-graphs $\phi \colon \sigma \to \tau$, it
canonically constructs a weighted stable $B$-graph $\pi$, together with
a contraction of $B$-graphs $\psi \colon \pi \to \rho$ and a combinatorial
morphism of weighted marked graphs $b \colon \pi \to \sigma$:
$$
\xymatrix{
B & \pi \ar[r]^\psi \ar[d]^b & \rho \ar[d]^a \\
A \ar[u]^\xi & \sigma \ar[r]^\phi & \tau
}
$$
We call $\pi$ the {\it stable pullback of $\rho$ under $\phi$}.
We will describe how to obtain $\pi$ from $\rho$, assuming that $\phi$ is
an elementary contraction (i.e. it contracts a single edge).

If $\phi$ contracts a loop adjacent to a vertex $v \in V_\tau$, we simply
reattach a loop at every preimage $v' \in a_V^{-1}(v)$ (and decrease its
genus by one). If $\phi$ contracts an edge $\{f, \bar f\}$ connecting
the vertices $v_1 = \partial_\sigma(f), v_2 = \partial_\sigma(\bar f)$,
let $v = \phi_V(v_1) = \phi_V(v_2)$ their common image in $\tau$, and let
$v' \in a_V^{-1}$ be any vertex in the preimage of $v$ in $\rho$. There
can be two cases:
\begin{enumerate}
\item Replace $v'$ by two vertices $v_1', v_2'$ connected by an edge
$\{f', \bar f'\}$; their class and genus are determined by the
corresponding vertex in $\sigma$:
$\alpha_\pi(v_i') = \xi(\alpha_\sigma(v_i))$ and
$g_\pi(v_i') = g_\sigma(v_i)$. A flag $f_1$ of $v$ is moved to $v_1'$ or
$v_2'$ according to its position in $\sigma$, i.e. according to whether
$\phi^F(a_F(f_1))$ is attached to $v_1$ or $v_2$; its weight remains
unchanged. Now if either $v_1'$ or $v_2'$ is unstable, we undo this
construction and skip to case (2). Otherwise, it remains to
define the maps: $\psi$ is the map contracting $\{f', \bar f'\}$;
the combinatorial morphism $b$ is given by sending $v_i'$ to $v_i$, and
by sending a flag $f_1 \neq f'$ of $v_i'$ to
$\left(\phi^F \circ a \circ (\psi^F)^{-1}\right)(f_1)$. Other than that,
$b$ agrees with $a$.
\item 
Assume in the above construction, the vertex $v_2'$ was unstable.
We leave $\rho$ unchanged, and let $b_V$ send $v'$ to $v_1$.
Let $f_1$ be a flag of $v'$; we set $b_F(f_1) =  \phi^F(a_V(f_1))$ if that
is a flag attached to $v_1$, otherwise $b_F(f_1) = f$, where $f$ defined
above is part of the edge connection $v_1$ and $v_2$.
\end{enumerate}
The same construction is iteratively applied to every such vertex $v$ to
obtain $\pi$.

Geometrically, the contractions $\phi$ corresponds to the inclusion of a
boundary component $\Mbar(\sigma)$ of the moduli space $\Mbar(\tau)$ associated
to $\tau$, and the stable pull-back constructs the boundary component of
$\Mbar(\rho)$ upon which the boundary component $\Mbar(\sigma)$ is naturally
mapped by morphism $\Mbar(\tau) \to \Mbar(\rho)$ associated to $a$.

\begin{PropDef}
Let $(\xi, a, \tau', \phi) \colon (A, \tau) \to (B, \sigma)$ and
$(\eta, b, \sigma', \psi) \colon (B, \sigma) \to (C, \rho)$ be morphisms
of weighted stable marked graphs. Then we define the composition
$(\eta, b, \sigma', \psi) \circ (\xi, a, \tau', \phi)
	\colon (A, \tau) \to (C, \rho)$
to be $(\eta \xi, ac, \tau'', \psi \xi)$ where $(c, \tau'', \xi)$ is
the stable pullback of $\sigma'$ under $\phi$.

This composition is associative, defining the category of weighted
stable marked graphs.
\end{PropDef}
\[
\xymatrix{
C & \tau'' \ar[r]^\xi \ar[d]^c & \sigma' \ar[d]^b \ar[r]^\psi & \rho \\
B \ar[u]^\eta & \tau' \ar[r]^\phi \ar[d]^a & \sigma \\
A \ar[u]^\xi & \tau
}
\]
We denote by $\GGG^w_s$ the category of weighted stable marked graphs,
and by $\AAA$ the category of semigroups with indecomposable zeros.

\subsection{Weighted stable maps indexed by graphs.}
As in \cite[section 3]{BehrendManin}, let $\VVV$ be the category of
smooth projective varieties over a field $k$. Consider the fibered product
$\VVV \GGG^w_s$ of categories
$$
\xymatrix{
\VVV \GGG^w_s \ar[r] \ar[d] & \GGG^w_s \ar[d] \\
\VVV \ar[r]^{H^+_2} & \AAA
}
$$
where $H^+_2$ is the functor that associates to $V$ the semigroup of
effective classes in $\CH^1(V)$. Objects of $\VVV \GGG^w_s$ are pairs
$(V, \tau)$ where $V$ is a smooth projective variety over $k$ and
$\tau$ is a weighted stable $H^+_2(V)$-graph.

For any weighted graph $\tau$ and any vertex $v \in V_\tau$, let
$F_v = \{f \in F_\tau | \partial_\tau(f) = v\}$ be the set of
flags attached to $v$, and $\AA_v = \AA|_{F_v}$ be their weight data.
\begin{Def}
A stable map of type $(V, \tau)$ for an object
$(V, \tau)$ in $\VVV \GGG^w_w$ is a collection of stable maps
$(C_v, x_v, f_v)$ to $V$ of type $(g(v), \AA_v, \alpha(v))$
for every $v \in V_\tau$, such that
$f_{\partial_\tau(i)}(x_i) = f_{\partial_\tau(j_\tau(i))}(x_{j_\tau(i)})$
for all flags $i$.
\end{Def}

For a scheme $T$ and $(V, \tau) \in \VVV \GGG^w_s$, let
$\Mbar(T)(V, \tau)$ be the groupoid of families of
weighted stable maps of type $(V, \tau)$ over $T$, and let
$\Mbar(T)$ be the groupoid of arbitrary weighted stable maps.

\begin{Thm}			\label{thm:graph-DMstacks}
For a fixed scheme $T$, $\Mbar(T)$ defines a 2-functor
\[
\Mbar(T)(\underline{\hphantom{W}}) \colon \VVV \GGG^w_s \to
\Mbar(T).
\]

For every base change $u \colon T' \to T$, the  pullback
$u^* \colon \Mbar(T) \to \Mbar(T')$ commutes with the functors
$\Mbar(T)(\underline{\hphantom{W}})$ and 
$\Mbar(T')(\underline{\hphantom{W}})$.

Finally, for fixed $(V, \tau, \alpha)$, the category of weighted
stable maps of type $(V, \tau, \alpha)$ is a proper algebraic Deligne-Mumford
stack $\Mbar(V, \tau, \alpha)$ of finite type.
\end{Thm}

Of course, the compatibility with base change in particular implies that
that $\Mbar(\Phi)$ for some morphism $\Phi$ in $\VVV\GGG^w_s$ induces
a morphisms between the stacks associated by $\Mbar$ to the source
and target; i.e. $\Mbar$ is a 2-functor from $\VVV\GGG^2_s$ to the
2-category of Deligne-Mumford stacks.\footnote{Implicitly, we passed
from the description of a stack as a category fibered in groupoids to
the description as a 2-functor to the 2-category of groupoids. See e.g.
\cite[Chapter V]{M} for a discussion of both viewpoints.}

The last claim of the theorem immediately follows from theorem
\ref{thm:DMstacks} and the fact that by definition it is a closed substack
of $\prod_{v \in V_\tau} \Mbar_{g(v), \AA(v)}(V, \alpha(v))$.

To prove the first and second claim of the
theorem, we need to prove the existence of a functorial push-forward
in $\Mbar(T)$ associated to every morphism 
$(\xi, a, \tau', \phi) \colon (V, \tau) \to (W, \sigma)$ in
$\VVV\GGG^w_s$, and show that they are compatible with base change.
Every morphism in $\VVV\GGG^w_s$ can be written as a composition
of elementary morphisms of one of the following types:
changing the target (I), increasing the weights (II), forgetting a tail
(III), combining tails (IV), complete combinatorial morphisms (V),
contracting an edge (VI) and contracting a loop (VII).
For complete combinatorial morphisms this is immediate (and
there is nothing to add to the discussion in \cite[section 3, case
IV]{BehrendManin}). All other cases
have already been treated in \ref{sect:elementary-morphisms} in the
case where the target is a one-vertex graph; the general case follows
immediately from this.

What is left to prove is that the associated morphism are compatible
with composition in the category of weighted stable marked graphs, i.e.
that it does not depend on the way we break up a morphism into a
composition of elementary morphisms.

For compositions of contractions with contractions, respectively of the
(inverses of) combinatorial morphisms with combinatorial morphisms
this is immediate, and the only interesting case to prove is the case
of the composition $(\xi, a)^{-1} \circ \phi$ of (the formal inverse of)
a combinatorial morphism
$(\xi, a) \colon (B, \rho) \to (A, \tau)$ and a contraction of $A$-graphs
$\phi \colon \sigma \to \tau$. In fact, the formation of stable
pull-back exactly makes sure that this compatibility holds, and the claim
follows easily by following every step of the stable pull-back construction.

\section{Graph-level description of virtual fundamental classes.}
					\label{sect:virtual-and-graphs}

To define Gromov-Witten invariants based on weighted stable maps, we need to
define virtual fundamental classes in the Chow ring
$A_*\left(\Mbar(V, \tau, \alpha)\right)$
of the moduli spaces. To formulate the required behavior with respect to
restriction to boundary components of the moduli space, 
we need to introduce the notion of isogenies of weighted stable graphs
and their cartesian isogeny diagrams. (We won't 
introduce the complete cartesian extended isogeny category as in 
\cite{BehrendManin}.)

\subsection{Isogenies of graphs.}
For our purposes, we need to refine the definition of an isogeny
as given in \cite[Definition 5.4]{BehrendManin}.

\begin{Def}
We say that the one-vertex $V$-graph $\sigma$ is a \emph{contraction of small
tails} of the one-vertex $V$-graph $\tau$ if it is obtained from $\tau$
by a sequence of steps, each forgetting a single tail, such that in
every step we are in
the situation of corollary \ref{cor:universal-curve} (the weight data
of $\tau$ is contained in a fine open chamber, and the weight of the
additional tail in $\sigma$ is small enough that changing it to zero
would not cross a wall of the fine chamber decomposition).
\end{Def}
This implies that the associated map $\Mbar(\tau) \to \Mbar(\sigma)$ is
flat, as it is a sequence of projection maps of the universal curve.

\begin{Def}
An isogeny $\Phi \colon \tau \to \sigma$ of weighted stable $A$-graphs 
is given by an injective map $\Phi^F \colon F_\sigma \to F_\tau$ of flags
and a surjective map $\Phi_V \colon V_\tau \to V_\sigma$ of vertices such
that the following conditions hold:
\begin{enumerate}
\item 
$\Phi^F$ commutes with the boundary maps $\partial_\tau, \partial_\sigma$,
i. e. for any flag $f \in F_\sigma$, we have
$\Phi_V(\partial_\tau(\Phi^F(f))) = \partial_\sigma(f)$.
\item 			\label{isog:genus}
For any vertex $v \in V_\sigma$, let $\tau_v$ be the subgraph of
$\tau$ that consists of all vertices send to $v$ by $\Phi_V$, and all edges
joining them. We require that
\begin{enumerate}
\item $g(v) = \sum_{w \in V_{\tau_v}} g(w) + \dim H^1(\abs{\tau_v})$ and
\item $\alpha(v) = \sum_{w \in V_{\tau_v}} \alpha(w)$
\end{enumerate}
\item
$\Phi^F$ respects the weights, i.e. $\AA_\tau \circ \Phi^F = \AA_\sigma$.
\item 					\label{isog:tail_forget}
For any $v \in V_\tau$, let $\tau_v$ be the one-vertex graph obtained from
$\tau$ by removing all other vertices, and cutting off the edges starting
from $v$ into a tail of weight 1;
let $\sigma_v$ be the graph obtained from $\tau_v$ by removing all
tails not in the image of $\Phi^F$. The condition is that $\sigma_v$ is
is a contraction of small tails of $\tau_v$.
\end{enumerate}
\end{Def}
In the geometric realizations of the graphs, an isogeny is given by
collapsing a set of disjoint closed connected subgraphs
$\abs{\tau_v} \subset \abs{\tau}$ consisting of edges and small tails
to a single vertex $v \in V_\sigma$.
It can be written as the composition of a morphism contracting small tails,
and a contraction as in definition \ref{def:contraction}.

\subsection{Cartesian isogeny diagrams}
Consider a stable $V$-graph $\sigma$ and its absolute stabilization
$a \colon \sigma^s \to \sigma$, as well as an isogeny of weighted modular
graphs $\Phi \colon \tau^s \to \sigma^s$. In \cite[section 5]{BehrendManin} the
pull-back $\tau = (\tau_i)_{i \in I}$ of $\sigma$ along $\Phi$ is constructed.
For each $i \in I$, the stable $V$-graph $\tau_i$ comes with a stabilization
morphism $a_i \colon \tau^s \to \tau_i$ and an isogeny
$\Phi_i \colon \tau_i \to \sigma$ such that the diagram
$$
\xymatrix{
\tau_i \ar[r]^{\Phi_i} 
& \sigma \\
\tau^s \ar[r]^\Phi \ar[u]^{a_i} & \sigma^s \ar[u]^b 
}
$$
commutes.

Its construction is as follows:\footnote{Unlike \cite[section
5]{BehrendManin}, we omit
the orbit map as well as the notion of an extended isogeny.}
To every edge $\{f, \bar f\}$ of $\sigma^s$ there is a long edge in
$\sigma$ consisting of edges
$\{f_1, \bar f_1\}, \dots, \{f_n, \bar f_n\}$ and vertices
$v_i = \partial_\sigma(\bar f_i) = \partial_\sigma(f_{i+1})$
such that
$b^F(f) = f_1$, $b^F(\bar f) = \bar f_n$ and the vertices $v_i$
are of genus 0 and have no further flags. We replace the edge
$\{\Phi^F(f), \Phi^F(\bar f)\}$ of $\tau^s$ by the same long
edge $\{f_1, \bar f_1\}, \dots, \{f_n, \bar f_n\}$.
Similarly, to every tail $f \in S_{\sigma^s}$ there is a long tail
$\{f_1, \bar f_1\}, \dots, \{f_n, \bar f_{n}\}$ of edges as above
and some number $k \ge 0$ of additional tails
$f_{n+1}, \dots, f_{n+k}$. The additional tails are attached to
the last vertex $v_n$ of the tail,
$\partial_\sigma(f_{n+i}) = v_n = \partial_\sigma(\bar f_n)$
for $1 \le i \le k$, and the
sum of weights is bounded as
$\sum_{1 \le i \le k} \AA(f_{n+i}) \le 1$.
Again we replace the tail $\Phi^F(f) \in S_{\sigma^s}$
with the same long tail, preserving the weights.

We thus obtain a weighted graph $\tau'$ with a combinatorial morphism
$a \colon \tau^s \to \tau'$ and an isogeny of graphs
$\Phi' \colon \tau' \to \sigma$.
Now let $I$ be the set of $V$-structures on $\tau'$ such that 
$\Phi'$ is an isogeny of weighted $V$-graphs. We get a set
$(\tau_i)_{i \in I}$ of $V$-graphs such that 
the induced morphism $a_i \colon \tau^s \to \tau_i$ is an absolute
stabilization, and $\Phi_i \colon \tau_i \to \sigma$ is an isogeny
of $V$-graphs.

The same construction can be made for a tuple $(\sigma_j)_{j \in J}$
of $V$-graphs with absolute stabilization 
morphisms $b_j \colon \sigma \to \sigma_j$. The formation of pull-back
then becomes compatible with composition.

\subsection{Expected properties.}
\begin{Def}
Let $\tau$ be a weighted stable $V$-graph, where $V$ is of pure
dimension $\dim V$, and has canonical class $\omega_V$.
We define the class $\beta(\tau)$, the Euler characteristic
$\chi(\tau)$, the genus $g(\tau)$ and the dimension $\dim(\tau)$ of
$\tau$ as
\begin{flalign*} 
\beta(\tau) &= \sum_{v \in V_\tau} \beta(v) \\
\chi(\tau) &= \chi(\abs{\tau}) - \sum_{v \in V_\tau} g(v) \\
g(\tau) &= 1 - \chi(\tau) \\
\dim (\tau) &= \chi(\tau)(\dim V - 3) - \beta(\tau)\cdot\omega_V
	       + \abs{S_\tau} - \abs{E_\tau}
\end{flalign*}
\end{Def}

We now fix $V$. An orientation will be a system of virtual fundamental
classes $J(V, \tau) \subset A_{\dim(V, \tau)}(\Mbar(V, \tau))$
for all stable $V$-graphs $\tau$ bounded by the characteristic, satisfying
the list of properties given below. 

\begin{itemize}
\item[(1)]
(\emph{Mapping to a point}). If $\tau$ is a graph of class zero,
and $\abs{\tau}$ is nonempty and connected, then
$$
J(V, \tau) = c_{g(\tau) \dim V} \left(R^1\pi_*f^* TV\right).
$$

\item[(2)] (\emph{Forgetting tails}).
Let $\Phi \colon \sigma \to \tau$ be a morphism of stable
$V$-graphs given by forgetting a small tail of $\sigma$, i.e.
such that $\tau$ is obtained from $\sigma$ by a contraction of a small
tail.
Then $M(\Phi)$ is flat, and we require
$$
J(V, \sigma) = \Mbar(\Phi)^* J(V, \tau).
$$

\item[(3)]
(\emph{Combining tails.}).
Let $\Phi \colon \sigma \to \tau$ be a morphism splitting
up a tail into several of them, i.~e. one that is induced by a
combinatorial morphism
$a \colon \tau \to \sigma$ combining several tails $f_1, \dots,
f_k \in S_\tau$ to a single tail $f \in S_\sigma$ with weight
$\AA_\sigma(f) = \sum_{i=1}^k \AA_\tau(f_i)$.
Then $\Mbar(\Phi)$ is a regular closed embedding, and the 
required condition is
$$
J(V, \sigma) = \Mbar(\Phi)^! J(V, \tau).
$$

\item[(4a)]
(\emph{Products}).
For any two stable $V$-graphs $\sigma, \tau$, let
$\sigma \times \tau$ be the disjoint union of the graphs
of $\sigma$ and $\tau$ with the obvious structure as a stable $V$-graph.
Then
$$
J(V, \sigma \times \tau) = J(V, \sigma) \times J(V, \tau).
$$

\item[(4b)]
 (\emph{Cutting edges}). Let $\Phi \colon \sigma \to \tau$
be a morphism obtained by cutting an edge $\{f, \bar f\}$ of $\sigma$ into
two tails. By abuse of notation, we identify the flags
$f, \bar f \subset F_\sigma$ with the corresponding tails
$f, \bar f \subset S_\tau$.
We obtain a cartesian square
\[
\xymatrix{
\Mbar(V, \sigma) \ar[d]^{\ev_f = \ev_{\bar f}} \ar[r]^{\Mbar(\Phi)}
& \Mbar(V, \tau) \ar[d]^{\ev_f \times \ev_{\bar f}} \\
V \ar[r]^\Delta & V \times V
}
\]
and require that 
$$J(V, \sigma) = \Delta^! J(V, \tau).$$

\item[(4c)] (\emph{Isogenies}).
Let $(\sigma_j)_{j \in J}$ be a tuple of $V$-graphs with absolute stabilization
$\sigma^s$ and $\tau^s \to \sigma^s$ an isogeny. Let
$(\tau_i)_{i \in I}$ be the tuple of $V$-graphs completing this to a
cartesian isogeny diagram. We obtain an induced commutative, but not
cartesian diagram
$$
\xymatrix{
\coprod_{i \in I} \Mbar(\tau_i) \ar[r] \ar[d] &
\coprod_{j \in J} \Mbar(\sigma_j) \ar[d] \\
\Mbar(\tau^s) \ar[r]^{\Mbar(\Phi)} &
\Mbar(\sigma^s)
}
$$
and thus an induced map
$$
h \colon \coprod_{i \in I} \Mbar(\tau_i) \to
\Mbar(\tau^s) \times_{\Mbar(\sigma^s)}
\coprod_{j \in J} \Mbar(\sigma_j).
$$
We require that 
$$
h_*\left(\sum_{i \in I} J(V, \tau_i)\right)
= \sum_{j \in J} \Mbar(\Phi)^! J(V, \sigma_j).
$$

\item[(5)] \emph{Kontsevich-stable maps.}
Assume that all weights satisfy $\AA(s) = 1$. Then 
$J(V, \tau)$ agrees with the definition of the
virtual fundamental class $J(V, \tilde \tau)$ for the underlying
stable $V$-graphs $\tilde \tau$ according to
\cite{Behrend:GW-alg, BF}.

\item[(6)] \emph{Reducing weights.}		
Let $\Phi \colon \sigma \to \tau$ be a morphism of
weighted stable
$V$-graphs obtained by reducing weights, i.~e. such that $\Phi$ is induced
by a combinatorial morphism $\tau \to \sigma$ that is the identity on the
modular graph structure, but such that $\AA_\sigma (f) \ge \AA_\tau (f)$
for all flags $f \in F_\tau = F_\sigma$. Then
$\Mbar(\Phi)$ is a reduction morphism, and we require that
$$
\Mbar(\Phi)_* \left(J(V, \sigma)\right) = \left(J(V, \tau)\right).
$$
\end{itemize}

\begin{Thm}			\label{thm:graph-level-virt-class}
There is a system of virtual fundamental classes satisfying all
properties listed in the previous section.
\end{Thm}

Note that (4a), (4b) and (4c) imply condition (4) of theorem
\ref{thm:virt_class}, whereas the other conditions for one-vertex graphs
are identical to the corresponding condition in \cite{BehrendManin}.

Of course, (1), (2) and (4a-c) are direct generalizations of properties
of the virtual fundamental classes in the non-weighted setting. The only
caveat is that for morphisms contracting or forgetting a tail, we
always have to assume the situation of corollary \ref{cor:universal-curve}.
This is to be expected: if we forget a tail of bigger weight, the forgetful
map factorizes via a non-trivial reduction morphism $\rho$. However, there
is no reason to assume that the virtual fundamental class is a pull-back of
a class via $\rho$.

As we already explained, we use (5) and (6) as the definition:
\begin{DefRem} 
For any weighted stable $V$-graph $\tau$, let $\tau^1$ be
the weighted stable $V$-graph obtained by setting all weights to 1,
let $w(\tau) \colon \tau \to \tau^1$ be the combinatorial morphism
increasing the weights, and $W(\tau) \colon \tau^1 \to \tau$ the induced
morphism in the category of weighted marked graphs. Then any
combinatorial morphism $\tau \to \sigma$ to a $V$-graph $\sigma$
with all weights equal to 1 factors uniquely via $w(\tau)$.
\end{DefRem}
By abuse of notation, we write $W(\tau) \colon \Mbar(V, \tau^1) \to
\Mbar(V, \tau)$ also for the induced map on moduli spaces, and
define $J(V, \tau)$ as
$$
J(V, \tau) := W(\tau)_* J(V, \tau^1)
$$
where the latter is as defined in \cite{Behrend:GW-alg, BF}.

We will now show how to obtain these properties from those listed in
Definition 7.1 in \cite{BehrendManin}, which have been verified for the
Behrend-Fantechi construction of the virtual fundamental class in
\cite{Behrend:GW-alg}. As a preparation, we need the following lemma:
\begin{Lem}				\label{lem:int-cartesian}
Let $\Phi \colon \sigma \to \tau$ be an isogeny of $V$-graphs, and let
$\Phi^1 \colon \sigma^1 \to \tau^1$ be the same morphism for the graphs
with weight 1.
Consider the commutative (but not necessarily cartesian) square
\[
\xymatrix{
\Mbar(V, \sigma^1) \ar[d]^{\Mbar(W(\sigma))} \ar[r]^{\Mbar(\Phi^1)}
& \Mbar(V, \tau^1) \ar[d]^{\Mbar(W(\tau))} \\
\Mbar(V, \sigma) \ar[r]_{\Mbar(\Phi)}
& \Mbar(V, \tau)
}
\]
and the induced morphism $h \colon \Mbar(V, \sigma^1) \to 
\Mbar(V, \sigma) \times_{\Mbar(V, \tau)} \Mbar(V, \tau^1)$.
Then $\Mbar(\Phi)^!$ and
$h_* \circ \Mbar(\Phi^1)^!$ yield the same orientation to the projection
\[
\Mbar(V, \sigma) \times_{\Mbar(V, \tau)} \Mbar(V, \tau^1)
	\to \Mbar(V, \tau^1).
\]
\end{Lem}
(By definition, an orientation of a morphism $f \colon X \to Y$ is an
element of the bivariant intersection theory $A^*(Y \to X)$, i.e.
in particular a morphism $A_*(X') \to A_*(Y')$ for every pull-back
$f' \colon X' \to Y'$ of $f$.)

We may assume that $\Phi$ is an elementary isogeny, so
we have one of the following two cases:
\begin{itemize}
\item \emph{Contraction of an edge}.
It is sufficient to consider the case where $\tau$ has only one vertex, so
both $\Mbar(\Phi)$ and $\Mbar(\Phi^1)$ are a gluing morphism as in 
proposition \ref{prop:gluing}. Consider the first case, where $\Phi$
contracts a non-looping edge (the other case follows similarly). An
object in the product consists of a pair of weighted stable maps
$((C_1, f_1), (C_2, f_2))$ of type $\sigma$ and $\tau^1$, respectively,
together with an isomorphism the reduction of $C_2$ to type $\tau$ with
the curve obtained by gluing the two components of $C_1$. Since the
sections cannot meet the node, this is only possible if $C_2$ already
consists of two components, which together form a weighted stable maps
of type $\sigma^1$. The induced map to $\Mbar(V, \sigma^1)$ is
an inverse to $h$, i.e. the above diagram is a cartesian square.

Both $\Mbar(\Phi)$ and $\Mbar(\Phi^1)$ are a codimension one regular
embedding with compatible normal bundle, and the claim follows by standard
intersection theory.

\item \emph{Contraction of a small tail}. 
In this case, both $\Mbar(\Phi)$ and $\Mbar(\Phi^1)$ are flat. The
orientation given by $\Mbar(\Phi)$ is the same as that of the projection
to the second factor of the product. Since $h$ is a
blow-up at a regularly embedded substack, we have $h_* \circ h^* = \id$,
and the assertion follows.
\end{itemize}

We proceed with the proof of theorem \ref{thm:graph-level-virt-class}.
\begin{itemize}
\item[(1)]
This follows from the same property
\cite[Definition 7.1, (1)]{BehrendManin} in the non-weighted case
and projection formula.

\item[(2)]
Consider the diagram of lemma \ref{lem:int-cartesian}:
\begin{align*}
\Mbar(\Phi)^! J(V, \tau)
& = \Mbar(\Phi)^! \Mbar(W(\tau))_* J(V, \tau^1)
& \text{(by definition)} \\
& = {p_1}_* \Mbar(\Phi)^! J(V, \tau^1)
& \text{(push-forward)} \\
& = {p_1}_* h_* \Mbar(\Phi^1)^! J(V, \tau^1)
& \text{(by lemma \ref{lem:int-cartesian})} \\
& = \Mbar(W(\sigma))_* J(V, \sigma^1)
& \text{(*)} \\
& = J(V, \sigma)
& \text{(by definition)}
\end{align*}
Here (*) holds by
\cite[Definition 7.1, (4)]{BehrendManin}.

\item[(4a)]
This is obvious from the same property for non-weighted graphs
\cite[Definition 7.1, (2)]{BehrendManin}.

\item[(4b)]
The natural map $\Mbar(V, \sigma^1) \to \Mbar(V, \tau^1)$
fits as an additional row on top of diagram given in condition
(4b), so
that all squares are cartesian.
Thus the claim follows from property \cite[Definition 7.1, (3)]{BehrendManin}
and push-forward.

\item[(4c)]
We may assume that $\abs{J} = 1$, so we are just dealing with a single
$V$-graph $\sigma$ and its absolute stabilization $\sigma^s$.

Consider
$\sigma^1$ and its absolute stabilization $(\sigma^1)^s$. By the universal
property of stabilization, the composition of the combinatorial morphisms
of weighted graphs $\sigma^s \to \sigma \to \sigma^1$ factors uniquely via
$(\sigma^1)^s$.
Similarly, for each $i \in I$ let $\tau_i^1$ be the corresponding
graphs with weights 1, and let, by some abuse of notation, $(\tau^1)^s$
be their common absolute stabilization; we obtain a combinatorial
morphism $\tau^s \to (\tau^1)^s$.

These morphisms can be completed to the following diagram of a cube:

\[
\xymatrix@R=12pt@C=7pt@M=1pt@W=5pt{
& {\coprod_i \tau_i} \ar[rrr]|{\coprod \Phi_i} \ar[ddd]
&&& \sigma \ar[ddd]
\\
{\coprod_i \tau_i^1} \ar[ru]^-{W(\tau_i)}
	\ar[rrr]|{\coprod \Phi_i^1} \ar[ddd]
&&& {\sigma^1} \ar[ru]|{W(\sigma)} \ar[ddd]
\\
\\
& {\tau^s} \ar[rrr]|\Phi &&& {\sigma^s}
\\
{(\tau^1)^s} \ar[ru] \ar[rrr]|{\Phi^1} &&& {(\sigma^1)^s} \ar[ru]
}
\]
More precisely, 
there exist unique contractions
$\Phi^1 \colon (\tau^1)^s \to (\sigma^1)^s$ and
$\Phi_i^1 \colon \tau_i^1 \to \sigma^1$ such that
\begin{itemize}
\item[(I)] 	
the top and bottom square are commutative in the category of
weighted marked graphs, and
\item[(II)]
the square in front is a cartesian isogeny diagram.
\end{itemize}

Assuming these claims, the desired property can be deduced from 
the corresponding property \cite[Definition 7.1, (5)]{BehrendManin}
by careful diagram computation:

Since none of the squares of the cube necessarily yield cartesian
squares of moduli spaces, we need to consider the products
$P_{\mathrm{back}}
	= \Mbar(\tau^s) \times_{\Mbar (\sigma^s)} \Mbar(V, \sigma)$, 
$P_{\mathrm{front}}
	= \Mbar((\tau^1)^s) \times_{\Mbar((\sigma^1)^s)} \Mbar(V, \sigma^1)$
and
$P_{\mathrm{diag}}
	= \Mbar(\tau^s) \times_{\Mbar(\sigma^s)} \Mbar(V, \sigma^1)$.
Let $h_{\mathrm{back}}$ and $h_{\mathrm{front}}$
be the induced map from the corresponding corner of the cube to
$P_{\mathrm{back}}$ and $P_{\mathrm{front}}$, respectively, and
$h_{\mathrm{d}\to\mathrm{b}}
\colon P_{\mathrm{diag}} \to P_{\mathrm{back}}$, 
$h_{\mathrm{f}\to\mathrm{b}}
\colon P_{\mathrm{front}} \to P_{\mathrm{back}}$ and
$h_{\mathrm{f}\to\mathrm{d}} \colon P_{\mathrm{front}} \to
P_{\mathrm{diag}}$ the maps induced by the commutative cube.
We obtain
\begin{align*}
\Mbar(\Phi)^! J(V, \sigma)
& = \Mbar(\Phi)^! \Mbar(W(\sigma))_* J(V, \sigma^1)
	& \text{(by definition)} \\
& = {h_{\mathrm{d}\to\mathrm{b}}}_* \Mbar(\Phi)^!  J(V, \sigma^1)
	& \text{(push-forward)} \\
& = {h_{\mathrm{d}\to\mathrm{b}}}_*
{h_{\mathrm{f}\to\mathrm{d}}}_* \Mbar(\Phi^1)^!  J(V, \sigma^1)
	& \text{(lemma \ref{lem:int-cartesian})} \\
& = {h_{\mathrm{f}\to\mathrm{b}}}_* {h_{\mathrm{front}}}_*
\sum_i J(V, \tau_i^1)
	& \text{(*)} \\
& = {h_\mathrm{back}}_* \sum_i W(\tau_i)_* J(V, \tau_i^1)
	& \\
& = {h_\mathrm{back}}_* \sum_i J(V, \tau_i),
	& \text{(by definition)}
\end{align*}
where (*) holds according to \cite[Definition 7.1, (5)]{BehrendManin}.
So it remains to prove the two claims above. 

The definition of $\Phi_i^1$ is obvious and necessarily unique, as the
graphs $\tau_i$ and $\tau_i^1$, as well as $\sigma_i$ and $\sigma_i^1$, 
are identical as marked graphs after forgetting the weights. Commutativity
of the top square is equivalent to the claim that the combinatorial morphism
$w(\tau_i) \colon \tau_i \to \tau_i^1$ is the stable pull-back
(see p.~\pageref{stable-pullback})
of $w(\sigma) \colon \sigma \to \sigma^1$ along $\Phi_i^1$, which
is equally obvious.

For the bottom square involving $\Phi^1$, we need to review the
construction of cartesian isogenies. Consider any tail $f \in S_{\sigma^s}$;
it corresponds to a long tail in $\sigma$ consisting of edges
$\{f_1, \bar f_1\}, \dots, \{f_n, \bar f_n\}$, of
vertices $v_1, \dots, v_n$ and of tails $f_{n+1}, \dots, f_{n+k}$ attached
to $v_n$. Its preimage $\Phi^F(f) \in S_{\tau^s}$ corresponds to an
identical long tail $\{\Phi_i^F(f_1), \Phi_i^F(\bar f_1)\}, ...$ etc. in
$\tau_i$. After adjusting the weights to one, we again see identical
long tails as part of $\sigma^1$ respectively $\tau_i^1$; these will have
identical stabilization in $(\sigma^1)^s$ resp. $(\tau^1)^s$. This
shows that $\Phi^1$ is uniquely determined on the stabilization of
this long tail.  The same discussion applies to any edge of $\sigma^s$
corresponding to a long edge in $\sigma^s$. 
Finally, any part of $\tau^s$ contracted by $\Phi$ will appear identically
in $\tau_i$, and thus in $\tau_i^1$ and $(\tau^1)^s$. Hence $\Phi^1$
will necessarily contract it, too.

We have thus constructed $\Phi^1$ so that the front square is a cartesian
isogeny diagram. At the same time, the above discussion
shows that the stable pull-back of $\sigma^s \to (\sigma^1)^s$ along
$\Phi^1$ will recover $\tau^s \to (\tau^1)^s$, i.e. the bottom
square is indeed commutative.

\item[(5)] This holds by definition.

\item[(6)] This follows from the definition and the fact that reduction
morphisms are compatible with composition (Proposition \ref{prop:reduction}).

\item[(3)]
By properties (4a) and (4b), we can consider only graphs
having a single vertex. Further, we may assume that the combinatorial
morphism $a$ combines exactly two tails $f_1, f_2 \in S_\tau$ to
a single tail $f = a_F(f_1) = a_F(f_2) \in S_\sigma$.

Let $\rho$ be the $V$-graph obtained from $\sigma^1$ by adding a second
vertex of class and genus zero, having two tails $f_1', f_2'$ of weight 1 and
one edge whose second flag connects it to the original vertex and replaces
the tail $f$; geometrically, we attach a tripod\footnote{a graph with a
single vertex and 3 tails} to the tail $f$.

The morphism $\rho \to \sigma^1$ induced by the combinatorial morphism
$\sigma^1 \to \rho$ gives an isomorphism of moduli spaces
$\Mbar(\rho) \to \Mbar(\sigma^1)$, which respects the virtual fundamental
classes by properties (1), (4a) and (4b).

There is a morphism $\Psi \colon \rho \to \tau^1$ contracting the edge
in $\rho$ and sending $f_i'$ to $f_i$. Thus we have the following
commutative diagram:
\[
\xymatrix{
\Mbar(\rho) \cong \Mbar(\sigma^1) \ar[r]^{\Mbar(\Psi)} \ar[d]^{W(\sigma)}
	& \Mbar(\tau^1) \ar[d]^{W(\tau)} \\
\Mbar(\sigma) \ar[r]^{\Mbar(\Phi)}
	& \Mbar(\tau)
}
\]
A discussion similar to the one in the proof of (4c) shows that this is
a cartesian square. Let $\Xi \colon \tau^1 \to \sigma^1$ be the morphism
obtained by forgetting the tail $f_1$ and mapping $f_2$ to $f$.
Then $\Mbar(\Psi)$ is a section of $\Mbar(\Xi)$, so
$\Mbar(\Psi)^! [\Mbar(\tau^1)]^\virt =
\Mbar(\Psi)^! \Mbar(\Xi)^* [\Mbar(\sigma^1)]^\virt =
[\Mbar(\sigma^1)]^\virt$. The desired equality follows by push-forward
and the vanishing of excess intersection.
\end{itemize}

\bibliography{../all}                      
\bibliographystyle{../alphaspecial}     

\end{document}